\documentclass[11pt]{article}

\usepackage{fullpage}

\usepackage{amsmath,amssymb,amsthm,url,amsfonts,mathrsfs,dsfont}

\usepackage{graphicx, epsfig}
\usepackage{verbatim}
\usepackage{subcaption}
\usepackage{ifthen}
\usepackage{multirow}
\usepackage{multicol}
\usepackage{float}
\usepackage{subeqnarray}
\usepackage{enumerate, enumitem}
\usepackage{cite}
\usepackage[usenames,dvipsnames]{color}
\usepackage[font=small]{caption}

\usepackage{pst-eucl}
\usepackage[locale=DE]{siunitx}

\usepackage[bookmarks=true,pageanchor,colorlinks,linkcolor=purple,anchorcolor=purple,citecolor=purple,urlcolor=purple]{hyperref}\makeatletter \g@addto@macro\UrlSpecials{\do\!{\newline}}
\usepackage{breakurl}

\usepackage{algorithm}
\usepackage{algpseudocode}

\algnewcommand{\Initialize}[1]{%
  \State \textbf{Initialize:}
  \Statex \hspace*{\algorithmicindent}\parbox[t]{.8\linewidth}{\raggedright #1}
}

\algdef{SE}[DOWHILE]{Do}{doWhile}{\algorithmicdo}[1]{\algorithmicwhile\ #1}%

\usepackage{makecell}

\usepackage{varwidth}

\def\hlinewd#1{%
  \noalign{\ifnum0=`}\fi\hrule \@height #1 \futurelet
   \reserved@a\@xhline}

\usepackage{ifthen}
\newboolean{showcomments}
\setboolean{showcomments}{true}

\newcommand{\ye}[1]{\ifthenelse{\boolean{showcomments}}
{ \textcolor{red}{(Ye says:  #1)}}{}}
\newcommand{\bose}[1]{  \ifthenelse{\boolean{showcomments}}
{ \textcolor{red}{(Bose says:  #1)}} {}  }
\newcommand{\lang}[1]{  \ifthenelse{\boolean{showcomments}}
{ \textcolor{red}{(Lang says:  #1)}} {}  }

\newcommand{\yp}{y^{\text{opt}}}
\newcommand{\Jp}{J^{\text{opt}}}

\newcommand{\beq}{\begin{equation}}
\newcommand{\eeq}{\end{equation}}
\newcommand{\beqa}{\begin{eqnarray}}
\newcommand{\eeqa}{\end{eqnarray}}
\newcommand{\beqan}{\begin{eqnarray*}}
\newcommand{\eeqan}{\end{eqnarray*}}

\renewcommand{\(}{\left (}
\renewcommand{\)}{\right )}
\renewcommand{\[}{\left [}
\renewcommand{\]}{\right]}

\newcommand{\vnorm}[1]{\left\|#1\right\|}

\newcommand\T{{\mathpalette\raiseT\intercal}}
\newcommand\raiseT[2]{\raisebox{0.25ex}{$#1#2$}\hspace{-0.1cm}}

\newcommand{\conv}{{\rm conv}}

\newcommand{\diag}{\mathop{\mathrm{diag}}}

\newcommand{\Rset}{\mathbb{R}}

\newcommand{\Acal}{{\cal A}}

\newcommand{\Dcal}{{\cal D}}

\newcommand{\Ical}{{\cal I}}

\newcommand{\Pcal}{{\cal P}}

\newcommand{\Vcal}{{\cal V}}

\newcommand{\Ycal}{{\cal Y}}

\newcommand{\Msf}{\sf{M}}

\newcommand{\Psf}{\sf{P}}

\newcommand{\argmin}{\mathop{\rm argmin}}

\newcommand{\bone}{\mathbf{1}}

\renewcommand{\v}[1]{{\mathbf{#1}}}
\newcommand{\ve}{\varepsilon}

\newcounter{l1}
\newcounter{l2}
\newcounter{l3}
\newcommand{\bdotlist}{\begin{list}{$\bullet$}{}}
\newcommand{\bboxlist}{\begin{list}{$\Box$}{}}
\newcommand{\bbboxlist}{\begin{list}{\raisebox{.005in}{{\tiny
$\blacksquare$ \ \ }}}{}}
\newcommand{\bdashlist}{\begin{list}{$-$}{} }
\newcommand{\blist}{\begin{list}{}{} }
\newcommand{\barablist}{\begin{list}{\arabic{l1}}{\usecounter{l1}}}
\newcommand{\balphlist}{\begin{list}{(\alph{l2})}{\usecounter{l2}}}
\newcommand{\bAlphlist}{\begin{list}{\Alph{l2}.}{\usecounter{l2}}}
\newcommand{\bdiamlist}{\begin{list}{$\diamond$}{}}
\newcommand{\bromalist}{\begin{list}{(\roman{l3})}{\usecounter{l3}}}

\newtheorem{theorem}{Theorem}
\newtheorem{lemma}{Lemma}

\newtheorem{remark}{Remark}

\renewcommand{\Rset}{\mathbf{R}}
\renewcommand{\bone}{\mathds{1}}

\usepackage{times}

\graphicspath{{Figures/}}

\begin{document}

\title{{On Robust Tie-line Scheduling \\ in Multi-Area Power Systems} \\ {\small Working paper}}
\author{{Ye Guo \qquad Subhonmesh Bose \qquad Lang Tong}
\thanks{\scriptsize
Y. Guo and L. Tong are with the School of Electrical and Computer Engineering, Cornell University, Ithaca, NY, USA. (Emails: {\{yg299,lt35\}@cornell.edu}). S. Bose is with the Dept. of Electrical and Computer Engineering at the University of Illinois at Urbana-Champaign, Urbana, IL, USA. (Email: {boses@illinois.edu}).}
}
\date{}

\maketitle

\begin{abstract}
The tie-line scheduling problem in a multi-area power system seeks to optimize tie-line power flows across areas that are independently operated by different system operators (SOs). In this paper, we leverage the theory of multi-parametric linear programming to propose algorithms for optimal tie-line scheduling respectively within a deterministic and a robust optimization framework. Aided by a coordinator, the proposed methods are proved to converge to the optimal schedule within a finite number of iterations. A key feature of the proposed algorithms, besides their finite step convergence, is that SOs do not reveal their dispatch cost structures, network constraints, or natures of uncertainty sets to the coordinator. The performance of the algorithms is evaluated using several power system examples.
\end{abstract}


\section{Introduction}
\label{sec:intro}


For historic and technical reasons, different parts of an interconnected power system and their associated assets are dispatched by different system operators (SOs). We call the geographical footprint within an SO's jurisdiction an \emph{area}, and transmission lines that interconnect two different areas as \emph{tie-lines}. Power flows over such tie-lines are generally scheduled 15 -- 75 minutes prior to power delivery. The report in \cite{MISO&PJM:10Overview} indicates that current scheduling techniques often lead to suboptimal tie-line power flows. The economic loss due to inefficient tie-line scheduling is estimated to the tune of \$73 million between the areas controlled by MISO and PJM alone in 2010. Tie-lines often have enough transfer capability to fulfill a significant portion of each area's power consumption \cite{White&Pike:11WP}. Thus they form important assets of multi-area power systems.

SOs from multiple areas typically cannot aggregate their dispatch cost structures and detailed network constraints to solve a \textit{joint} optimal power flow problem. Therefore, distributed algorithms have been proposed. Prominent examples include \cite{Kim&Baldick:97TPS, ConejoAguado98TPS, Bakirtzis&Biskas:03TPS} that adopt the so-called \emph{dual decomposition} approach. These methods are iterative, wherein each SO optimizes the grid assets within its area, given the Lagrange multipliers associated with inter-area constraints. Typically, a coordinator mediates among the SOs and iteratively updates the multipliers. Alternative \emph{primal decomposition} approaches are also proposed in \cite{GuoTongetc:CRP_TPWRS, Li&Wu&Zhang&Wang::15TPS, Zhao&LitvinovZheng:14TPS}. Therein, the primal variables of the optimization problem are iteratively updated, sometimes requiring the SO of one area to reveal part of its cost structure and constraints to the SO of another area or a coordinator.

%

Traditionally, solution techniques for the tie-line scheduling problem assume that the SOs and/or the coordinator has perfect knowledge of the future demand and supply conditions at the time of scheduling. Such assumptions are being increasingly challenged with the rapid adoption of distributed energy resources in the distribution grid and variable renewable generation like wind and solar energy in the bulk power systems. Said differently, one must explicitly account for the uncertainty in demand and supply in the tie-line scheduling problem. To that end, \cite{ahmadi2014multi,JiZhengTong:16TPS} propose to minimize the expected aggregate dispatch cost and \cite{LiWu_RobustMAED2016} propose to minimize the maximum of that cost. In this paper, we adopt the latter paradigm -- the \emph{robust} approach.

\subsection*{Our contribution}
With the system model in Section \ref{sec:systemModel}, we first formulate the deterministic tie-line scheduling problem in Section \ref{sec:deterministic}, where we propose an algorithm to solve this deterministic problem that draws from the theory of multiparametric programming \cite{borrelli2003constrained}. The key feature of our algorithm is that a coordinator can produce the optimal tie-line schedule upon communicating only \emph{finitely many times} with the SO in each area. In contrast to \cite{Zhao&LitvinovZheng:14TPS}, our method does not require SOs to reveal their cost structures nor their constraints to other SOs or to the coordinator. In Section \ref{sec:robust}, we formulate the robust counterpart of the tie-line scheduling problem. We then propose a technique that alternately uses the algorithm for the deterministic variant and a mixed-integer linear program to solve the robust problem. Again, our technique is proved to converge to the optimal robust tie-line schedule that requires the coordinator to communicate finitely many times with each SO. Also, SOs are not required to reveal the nature and range of the values the uncertain demand and available supply can take. Our proposed framework thus circumvents the substantial communication burden of the method proposed in\cite{LiWu_RobustMAED2016} towards the same problem. We remark that \cite{LiWu_RobustMAED2016} adopts the \emph{column-and-constraint generation} technique described in \cite{zeng2013solving} that requires SOs to reveal part of their network constraints, costs and ranges of demand and available renewable supply to the coordinator. We empirically demonstrate the performance of our algorithm in Section \ref{sec:simulation} and conclude in Section \ref{sec:conclusion}.

%



\section{System model}
\label{sec:systemModel}

\begin{figure}[h]
  \centering
  \vspace{0.2cm}
  \includegraphics[width=0.57\textwidth]{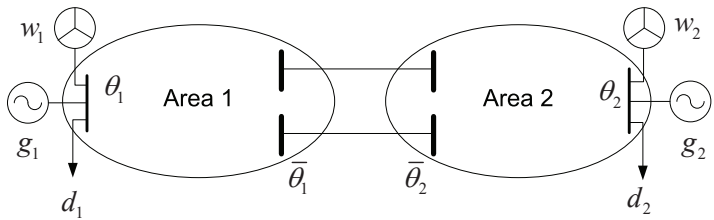}
  \caption{\small An illustration of a two-area power system.}\label{fig:1}
\end{figure}

To formulate the tie-line scheduling problem, we begin by describing the model for multi-area power systems. Throughout, we restrict ourselves to a two-area power system, pictorially represented in Figure \ref{fig:1} for the ease of exposition. The model and the proposed methods can be generalized for tie-line scheduling among more than two areas.

For the power network in each area, we distinguish between two types of buses: the \emph{internal} buses and the \emph{boundary} buses. The boundary ones in each area are connected to their counterparts in the other area via tie-lines. Internal buses do not share a connection to other areas. Assume that each internal bus has a dispatchable generator, a renewable generator, and a controllable load\footnote{While we assume that all loads are controllable, uncontrollable load at any node can be easily modeled by letting the limits on the allowable power demand at that node to be equal.}. Boundary buses do not have any asset that can inject or extract power. Such assumptions are not limiting in that one can derive an equivalent power network in each area that adheres to these assumptions.

Let the power network in area $i$ be comprised of $n_i$ internal buses and $\overline{n}_i$ boundary buses for each $i=1,2$. We adopt a linear DC power flow model in this paper.\footnote{See \cite{
Erseghe:15TPS, Magnússon2015ADMMACOPF}, and the references therein for solution approaches for a multi-area ACOPF problem.}
 This approximate model sets all voltage magnitudes to their nominal values, ignores transmission line resistances and shunt reactances, and deems differences among the voltage phase angles across each transmission line to be small.
Consequently, the real power injections into the network is a linear map of voltage phase angles (expressed in radians) across the network.
To arrive at a mathematical description, denote by $g_i \in \Rset^{n_i}$, $w_i \in \Rset^{n_i}$, and $d_i \!\in\! \Rset^{n_i}$ as the vectors of (real) power generations from dispatchable generators, renewable generators, and controllable loads, respectively. Let $\theta_i \in \Rset^{n_i}$ and $\overline{\theta}_i \in \Rset^{\overline{n}_i}$ be the vectors of voltage phase angles at internal and boundary buses, respectively. Then, the power flow equations are given by
\begin{align}
\begin{pmatrix}
\v{B}_{11} & \v{B}_{1\bar{1}}& & \\
\v{B}_{\bar{1}1} & \v{B}_{\bar{1}\bar{1}} & \v{B}_{\bar{1}\bar{2}} & \\
 & \v{B}_{\bar{2}\bar{1}} & \v{B}_{\bar{2}\bar{2}} & \v{B}_{\bar{2}2} \\
& & \v{B}_{2\bar{2}} & \v{B}_{22}
\end{pmatrix}
\begin{pmatrix}
\theta_{1}\\
\overline{\theta}_{1}\\
\overline{\theta}_{2}\\
\theta_{2}
\end{pmatrix}
=
\begin{pmatrix}
g_{1}+w_1-d_{1}\\
0\\
0\\
g_{2}+w_2-d_{2}
\end{pmatrix}.
\label{eq:ged_dclf}
\end{align}
Non-zero entries of the coefficient matrix depend on reciprocals of transmission line reactances, the unspecified blocks in that matrix are zeros. Throughout, assume that one of the boundary buses in area 1 is set as the slack bus for the two-area power system. That is, the voltage phase angle at said bus is assumed zero.

Power injections from the supply and demand assets at the internal buses of area $i$ are constrained as
\begin{align}
\underline{G}_i \leq g_i \leq \overline{G}_i,  \ \ 0 \leq w_i \leq \overline{W}_i, \ \ \underline{D}_i \leq d_i \leq \overline{D}_i.
\label{eq:gdw.const}
\end{align}
The inequalities are interpreted elementwise. The lower and upper limits on dispatchable generation $\underline{G}_i, \overline{G}_i$ are assumed to be known at the time when tie-line flows are being scheduled. Our assumptions on the available renewable generation $\overline{W}_i$ and the limits on the demands $[\underline{D}_i, \overline{D}_i]$ will vary in the subsequent sections. In Section \ref{sec:deterministic}, we assume that these limits are known and provide a distributed algorithm to solve the deterministic tie-line scheduling problem. In Section \ref{sec:robust}, we formulate the robust counterpart, where these limits are deemed uncertain and vary over a known set. We then describe a distributed algorithm to solve the robust counterpart.

The power transfer capabilities of transmission lines within area $i$ are succinctly represented as
\begin{align}
\v{H}_i \theta_i + \v{\overline{H}}_i \overline{\theta}_i \leq f_i
\label{eq:defHi}
\end{align}
for each $i=1,2$. Here, $\v{H}_i$ and $\v{\overline{H}}_i$ define the branch-bus admittance matrices, and $f_i$ models the respective transmission line capacities. Similarly, the transfer capabilities of tie-lines joining the two areas assume the form
\begin{align}
\v{\overline{H}}_{12} \overline{\theta}_1 + \v{\overline{H}}_{21} \overline{\theta}_2 \leq f_{12}.
\label{eq:defH12}
\end{align}
Again, $\v{\overline{H}}_{12}$, $\v{\overline{H}}_{21}$ denote the relevant branch-bus admittance matrices and $f_{12}$ models the tie-line capacities.

Finally, we describe the cost model for our two-area power system. For respectively procuring $g_i$ and $w_i$ from dispatchable and renewable generators, and meeting a demand of $d_i$ from controllable loads, let the dispatch cost in area $i$ be given by
\begin{align}
\[P^g_i\]^\T g_i + \[P^w_i\]^\T \( \overline{W}_i - w_i \) + \[P^d_i\]^\T \( \overline{D}_i - d_i \).
\label{eq:defP}
\end{align}
We use the notation $v^\T$ to denote the transpose of any vector or matrix $v$.
The linear cost structure in the above equation is reminiscent of electricity market practices in many parts of the U.S. today. The second summand models any spillage costs associated with renewable generators. The third models the disutility of not satisfying all demands.



\section{The deterministic tie-line scheduling problem}
\label{sec:deterministic}

Tie-line flows are typically scheduled ahead of the time of power delivery. The lead time makes the supply and demand conditions uncertain during the scheduling process. Within the framework of our model, the available capacity in renewable supply and lower and upper bounds on power demands, i.e., $\overline{W}_i, \underline{D}_i, \overline{D}_i$, can be uncertain. In this section, we ignore such uncertainty and formulate the deterministic tie-line scheduling problem, wherein we assume perfect knowledge of $\overline{W}_i$, $\underline{D}_i$ and $\overline{D}_i$ to decide the dispatch in each area and the tie-line flows. Our discussion of the deterministic version will serve as a prelude to its robust counterpart in Section \ref{sec:robust}.

To simplify exposition, consider the following notation.
\begin{align*}
x_i := \(g_{i},w_i,d_i,\theta_i \)^\T,  \ \ \xi_i  := \(\overline{W}_i, \underline{D}_i, \overline{D}_i \)^\T, \ \ y  := \( \overline{\theta}_1, \overline{\theta}_2 \)^\T
\end{align*}
for $i = 1,2$. The above notation allows us to succinctly represent the constraints \eqref{eq:ged_dclf} -- \eqref{eq:defHi} as
\begin{align*}
\v{A}^x_i x_i+ \v{A}^\xi_i \xi_i + \v{A}^y_i y \leq b_i
\end{align*}
for each $i = 1,2$ and suitably defined matrices $\v{A}^x_i, \v{A}^\xi_i, \v{A}^y_i$ and vector $b_i$. Denote by $m_i$ the number of inequality constraints in the above equation. Next, we describe transmission constraints on tie-line power flows in \eqref{eq:defH12} as
\begin{align*}
y \in \Ycal \subset \Rset^{Y}.
\end{align*}
Without loss of generality, one can restrict $\Ycal$ to be a polytope\footnote{Assuming the power network to be connected, the modulus of the phase angle of any bus can be constrained to lie within the sum of admittance-weighted transmission line capacities connecting that bus to the slack bus.}.
Finally, the cost of dispatch in area $i$, as described in \eqref{eq:defP}, can be written as
$$ c_i(x_i, \xi_i) := c^0_i + [c^x_i]^\T \ x_i + [ c^\xi_i ]^\T \ \xi_i$$
for scalar $c^0_i$ and vectors $c^x_i$, $c^\xi_i$.

Equipped with the above notation, we define the deterministic tie-line scheduling problem as follows.
\begin{equation}
\begin{alignedat}{8}  & \underset{x_1, x_2, y}{\text{minimize}}
	 & &  \left[ c_1\(x_1, \xi_1\) + c_2\(x_2, \xi_2\) \right], \\
 & \text{subject to}   \quad
	&& \v{A}^x_i x_i + \v{A}^\xi_i \xi_i + \v{A}^y_i y \leq b_i,  \ i = 1,2,\\
	&&& y \in \Ycal.
 \end{alignedat} \label{eq:detProb}
\end{equation}
%


\subsection{Distributed solution via critical region exploration}

The structure of the optimization problem in \eqref{eq:detProb} lends itself to a distributed solution architecture that we describe below. Our proposed technique is similar in spirit to the critical region projection method described in \cite{GuoTongetc:CRP_TPWRS}.\footnote{The cost structure in \cite{GuoTongetc:CRP_TPWRS} is quadratic; the linear cost case does not directly follow from \cite{GuoTongetc:CRP_TPWRS}.} We assume that each area is managed by a system operator (SO), and a \emph{coordinator} mediates between the SOs. Assume that the SO of area $i$ (call it SO$_i$) knows the dispatch cost $c_i$ and the linear constraint involving $x_i, \xi_i, y$ in \eqref{eq:detProb} in area $i$, and that SOs and the coordinator all know $\Ycal$.

Our algorithm relies on the properties of \eqref{eq:detProb} that we describe next. To that end, notice that \eqref{eq:detProb} can be written as
\begin{align}
\underset{y \in \Ycal}{\text{minimize}} \ {J^*\(y, \xi_1, \xi_2\)} := {J_1^*\(y, \xi_1\) + J_2^*\(y, \xi_2\)},
\label{eq:detProb.2}
\end{align}
where
\begin{equation} \label{eq:areaProb}
J_i^*\(y, \xi_i\) := \underset{x_i}{\text{minimum}}
\quad c_i\(x_i, \xi_i\), \quad \text{subject to}
\quad \v{A}^x_i x_i + \v{A}^\xi_i \xi_i + \v{A}^y_i y \leq b_i.
\end{equation}
Assume throughout that all optimization problems parameterized by $y$ is feasible for each $y \in \Ycal$. Techniques from \cite{LiWu_RobustMAED2016} can be leveraged to shrink $\Ycal$ appropriately, otherwise.
The optimization problem in \eqref{eq:areaProb} is a multi-parametric linear program, linearly parameterized in $\(y, \xi_i\)$ on the right-hand side\footnote{The problem in \eqref{eq:areaProb} reformulated using the so-called epigraph form yields a multi-parametric program that is classically recognized as one linearly parameterized on the right-hand side.}. Such optimization problems are well-studied in the literature. For example, see \cite{borrelli2003constrained}. Relevant to our algorithm is the structure of the parametric optimal cost $J_i^*$. Describing that structure requires an additional notation. We say that a finite collection of polytopes $\{ \Pcal^1, \ldots, \Pcal^\ell \}$ define a \emph{polyhedral partition} of $\Ycal$, if no two polytopes intersect except at their boundaries, and their union equals $\Ycal$. With this notation, we now record the properties of $J_i^*$ in the following lemma.
\begin{lemma}
\label{lemma:pAffinePoly}
$J_i^*(y, \xi_i)$ is piecewise affine and convex in $y \in \Ycal$. Sets over which $J_i^*( \cdot, \xi_i)$ is affine define a polyhedral partition of $\Ycal$.
\end{lemma}
The proof is immediate from \cite[Theorem 7.5]{borrelli2003constrained}. Details are omitted for brevity. We refer to the polytopes in the polyhedral partition of $\Ycal$ induced by $J_i^*(\cdot, \xi_i)$ as \emph{critical regions}. Recall that the feasible set of \eqref{eq:areaProb} is described by a collection of linear inequalities. Essentially, each critical region corresponds to the subset of $\Ycal$ over which a specific set of these inequality constraints are active -- \textit{i.e.}, are met with equalities -- at an optimal solution of \eqref{eq:areaProb}.

A direct consequence of the above lemma is that the \emph{aggregate cost} $J^*(\cdot, \xi_1, \xi_2)$ is also piecewise-affine and convex. Sets over which this cost is affine define a polyhedral partition of $\Ycal$. The polytopes of that partition -- the critical regions -- are precisely the non-empty intersections between the critical regions induced by $J_1^*(\cdot, \xi_1)$ and those by $J_2^*(\cdot, \xi_2)$. The relationship between the critical regions induced by the various piecewise affine functions are illustrated in Figure \ref{fig:criticalRegion}.
In what follows, we develop an algorithm wherein the coordinator defines a sequence of points in $\Ycal$ towards optimizing the aggregate cost. In each step, it relies on the SOs to identify their respective critical regions and the affine descriptions of their optimal costs at these iterates. That is, SO$_i$ can compute the critical region $\Pcal^y_i$ that contains $y \in \Ycal$ and the affine description $\[\alpha^y_i\]^\T z + \beta^y_i$ of its optimal dispatch cost $J_i^*\(z, \xi_i\)$ over $z \in \Pcal^y_i$ by parameterizing the linear program described in \eqref{eq:areaProb}\footnote{The critical region containing $y\in\Ycal$ is unique, except when $y$ lies at the boundary of critical regions. In that event, assume that the SO returns one of the critical regions containing $y$.}. We relegate the details of this step to Appendix \ref{sec:CRAffine} to maintain continuity of presentation. For any $y \in \Ycal$, we assume in the sequel that the coordinator can collect this information from the SOs to construct the critical region $\Pcal^y$ induced by the aggregate cost containing $y$ and its affine description $\[\alpha^y\]^\T z + \beta^y$ for $z\in \Pcal^y$, where
\begin{align}
\begin{aligned}
\Pcal^y := \Pcal_1^y \cap \Pcal_2^y, \ \
\alpha^y := \alpha_1^y  + \alpha_2^y, \ \
\beta^y := \beta_1^y + \beta_2^y.
\label{eq:alphaBeta}
\end{aligned}
\end{align}
%

\definecolor{newColorLine}{gray}{0.9}
\definecolor{darkYellowGreen}{RGB}{123, 164, 40}
\definecolor{darkOliveGreen}{rgb}{0.33, 0.42, 0.18}
\definecolor{darkTerracotta}{rgb}{0.8, 0.31, 0.36}
\definecolor{deepCarminePink}{rgb}{0.84, 0.09, 0.41}
\definecolor{myRed}{rgb}{0.84, 0.09, 0.41}

\psset{unit=0.23in}

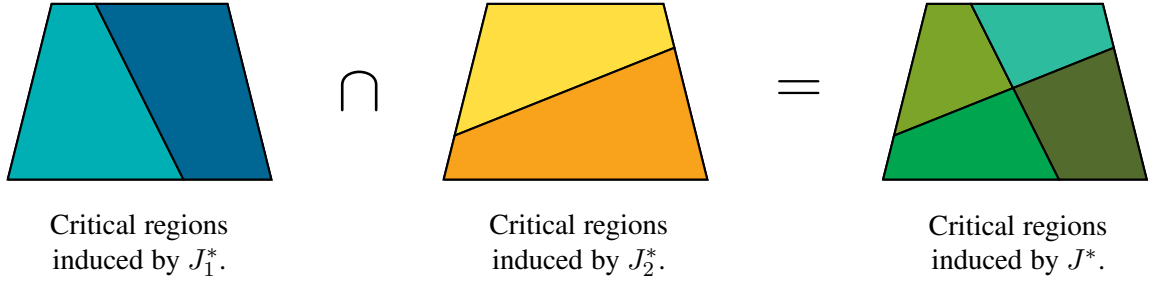
\begin{figure}[t!]
\centering
\begin{pspicture}[showgrid=false](-1, -3)(27, 4)

    \pstGeonode[PointName=none, PointSymbol=none]
        (0,0){A1}
        (1,4){A2}
        (5,4){A3}
        (6,0){A4}
        (2,4){B1}
        (4,0){B2}
	(0.25,1){B3}
	(5.25,3){B4}
    \pstInterLL[PointName=none,PointSymbol=none]{B1}{B2}{B3}{B4}{C}

    \rput(10,0){
        \pstGeonode[PointName=none, PointSymbol=none]
            (0,0){SA1}
            (1,4){SA2}
            (5,4){SA3}
            (6,0){SA4}
            (2,4){SB1}
            (4,0){SB2}
            (0.25,1){SB3}
            (5.25,3){SB4}
    \pstInterLL[PointName=none,PointSymbol=none]{SB1}{SB2}{SB3}{SB4}{SC}
    }

    \rput(20,0){
        \pstGeonode[PointName=none, PointSymbol=none]
            (0,0){TA1}
            (1,4){TA2}
            (5,4){TA3}
            (6,0){TA4}
            (2,4){TB1}
            (4,0){TB2}
            (0.25,1){TB3}
            (5.25,3){TB4}
    \pstInterLL[PointName=none,PointSymbol=none]{TB1}{TB2}{TB3}{TB4}{TC}
    }

    \psline (A1)(A2)(A3)(A4)(A1)
    \rput(3, -1.5){\shortstack[c]{Critical regions \\ induced by $J_1^*$.}}

    \rput(8, 2){\Huge $\cap$}

    \psline (SA1)(SA2)(SA3)(SA4)(SA1)
    \rput(13, -1.5){\shortstack[c]{Critical regions \\ induced by $J_2^*$.}}

    \rput(18, 2){\Huge $=$}

    \psline (TA1)(TA2)(TA3)(TA4)(TA1)
    \rput(23, -1.5){\shortstack[c]{Critical regions \\ induced by $J^*$.}}

    \pspolygon[fillstyle=solid, fillcolor=TealBlue](A1)(A2)(B1)(B2)
    \pspolygon[fillstyle=solid, fillcolor=MidnightBlue](B1)(A3)(A4)(B2)

    \pspolygon[fillstyle=solid, fillcolor=Goldenrod](SB3)(SA2)(SA3)(SB4)
    \pspolygon[fillstyle=solid, fillcolor=YellowOrange](SB3)(SB4)(SA4)(SA1)

    \pspolygon[fillstyle=solid, fillcolor=darkYellowGreen](TA2)(TB1)(TC)(TB3)
    \pspolygon[fillstyle=solid, fillcolor=SeaGreen](TB1)(TA3)(TB4)(TC)
    \pspolygon[fillstyle=solid, fillcolor=darkOliveGreen](TC)(TB4)(TA4)(TB2)
    \pspolygon[fillstyle=solid, fillcolor=Green](TB3)(TC)(TB2)(TA1)

\end{pspicture}
\caption{A pictorial representation of the critical regions induced by the areawise parametric optimal costs $J_1^*(\cdot, \xi_1), J_2^*(\cdot, \xi_2)$, and the aggregate cost $J^*(\cdot, \xi_1, \xi_2)$. The trapezoids represent $\Ycal$. Differently shaded polytopes indicate different critical regions.} \label{fig:criticalRegion}
\end{figure}

In presenting the algorithm, we assume that the coordinator can identify the \emph{lexicographically smallest} optimal solution of a linear program. A vector $a$ is said to be lexicographically smaller than $b$, if at the first index where they differ, the entry in $a$ is less than that in $b$. See \cite{dantzig2016linear} for details on such linear programming solvers. When a linear program does not have a unique optimizer\footnote{A linear program has non-unique optimizers when it is \emph{dual degenerate}. See \cite{dantzig2016linear} for details.}, such a choice provides a tie-breaking rule.
The final piece required to state and analyze the algorithm is an optimality condition that is both necessary and sufficient for a candidate minimizer of \eqref{eq:detProb.2}. Stated geometrically, $y^* \in \Ycal$ is a minimizer of \eqref{eq:detProb.2} if and only if
\begin{align}
0 \in \partial J^*( y^*, \xi_1, \xi_2) + N_\Ycal(y^*).
\label{eq:optCriterion}
\end{align}
The first set on the right-hand side of \eqref{eq:optCriterion} is the sub-differential set of the aggregate cost $J^*(\cdot, \xi_1, \xi_2)$ evaluated at $y^*$ \footnote{We use the sub-differential characterization as opposed to the familiar gradient condition for optimality since $J^*(\cdot, \xi_1, \xi_2)$ is piecewise affine and may not be differentiable everywhere in $\Ycal$.}. And, the second set denotes the normal cone to $\Ycal$ at $y^*$. The addition stands for a set-sum.

Algorithm \ref{alg:CRP} delineates the steps for the coordinator to solve the deterministic tie-line scheduling problem. In our algorithm, $\| v^* \|_2$ denotes the Euclidean norm of $v^*$. If $\Dcal := \{ \alpha^1, \ldots, \alpha^{\ell_D}\}$ and $N_\Ycal\(y^*\) := \{ z \ | \ \v{K}^y z \geq 0 \}$, then computing the least-square solution $v^*$ amounts to solving the following convex quadratic program.
\begin{eqnarray}
{\text{minimize}} \quad \frac{1}{2}\vnorm{v}_2^2,
\quad \text{subject to}
\quad v = \sum_{j=1}^{\ell_D} \eta_j \alpha^j + \zeta, \ \  \bone^\T \ \eta = 1, \ \ \eta \geq 0, \ \ \v{K}^y \zeta \geq 0
\label{eq:ifopt}
\end{eqnarray}
%
over the variables $v \in \Rset^{\overline{n}_1 + \overline{n}_2}$, $\eta \in \Rset^{\ell_D}$, and $\zeta\in\Rset^{\ell_N}$, where $\bone$ is a vector of all ones, and $\v{K}^y \in \Rset^{\( \overline{n}_1 + \overline{n}_2\) \times \ell_N}$.
\begin{algorithm}
	\caption{Solving the deterministic tie-line scheduling problem.}
	\label{alg:CRP}
	\begin{algorithmic}[1]
		\Initialize{$y \gets $ any point in $\Ycal$, $J^* \gets \infty$, \\$\Dcal \gets $ empty set, $\ve \gets$ small positive number.}
		\Do
	        		\State Communicate with the SOs to obtain $\Pcal^y$ and $\alpha^y, \beta^y$.
			\State Minimize $\[\alpha^y\]^\T z + \[\beta^y\]$ over $\Pcal^y$. \label{step:min}
			\State $\yp \gets$ lexicographically smallest minimizer in step \ref{step:min}.
			\State $\Jp \gets$ optimal cost in step \ref{step:min}.
        			\If {$ \Jp < J^*$,}
				\State $y^* \gets \yp$, $J^* \gets \Jp$, $\Dcal \gets \{ \alpha^y \}$.
        			\Else
        				\State $\Dcal \gets \Dcal \cup \{ \alpha^y \}$.

			\EndIf
				\State $v^* \gets \argmin_{v \in \conv(\Dcal) + N_{\Ycal}(y^*)} \| v \|_2^2$.
				\State $y \gets \yp - \ve v^*$. \label{step:yUpdate}

		\doWhile{$v^* \neq 0$.}
	\end{algorithmic}
\end{algorithm}
%


\subsection{Analysis of the algorithm}

The following result characterizes the convergence of Algorithm \ref{alg:CRP}. See Appendix \ref{sec:finiteTime} for its proof.
\begin{theorem}
\label{thm:finiteTime}
Algorithm \ref{alg:CRP} terminates after finitely many steps, and $y^*$ at termination optimally solves \eqref{eq:detProb.2}.
\end{theorem}
The above result fundamentally relies on the fact that each time the variable $y$ is updated, it belongs to a critical region (induced by the aggregate cost) that the algorithm has not encountered so far. And, there are only finitely many such critical regions. That ensures termination in finitely many steps. Each time the algorithm ventures into a new critical region, we store the optimizer and the optimal cost over that critical region in the variables $\yp$ and $\Jp$. Forcing the linear program to choose the lexicographically smallest optimizer always picks a unique vertex of the critical region as $\yp$. Unless $\Jp$ improves upon the cost at $y^*$, we ignore the new point $\yp$. However, the exploration of the new critical region provides a possibly new sub-gradient of the aggregate cost at $y^*$. The sub-differential set at $y^*$ is given by the convex hull of the sub-gradients of the aggregate cost over all critical regions that $y^*$ is a part of.
The set $\Dcal$ we maintain is such that $\conv(\Dcal)$ is a \emph{partial} sub-differential set of the aggregate cost at $y^*$. Notice that
$$ \conv(\Dcal) \subseteq \partial J^*(y^*, \xi_1, \xi_2)$$
throughout the algorithm. Therefore, any $y^*$ that meets the termination criterion of the algorithm automatically satisfies \eqref{eq:optCriterion}. As a result, such a $y^*$ is an optimizer of \eqref{eq:detProb.2}.

The proposed technique is attractive in that each SO only needs to communicate finitely many times with the coordinator for the latter to reach an optimal tie-line schedule. Further, each SO$_i$ can compute its optimal dispatch $x_i^*$ by solving \eqref{eq:areaProb} with $y^*$. A closer look at the nature of the communication between the SOs and the coordinator reveals that an SO will not have to disclose the complete cost structure nor a complete description of the constraints within its area to the coordinator.

\begin{remark}
\label{rem:pieceLine}
Algorithm \ref{alg:CRP} allows the coordinator to minimize
$$F(y) := F_1(y) + F_2(y)$$
in a distributed manner, where $F_i: \Ycal \to \Rset$ satisfies two properties. First, it is piecewise affine and convex. Second, given any $y \in \Ycal$, SO$_i$ can compute an affine segment containing that $y$. While we do not explicitly characterize how fast the algorithm converges to its optimum, one can expect the number of steps to convergence to grow with the number of critical regions so induced. However, we do not expect our algorithm to explore all such critical regions on its convergence path.
\end{remark}


\subsection{A pictorial illustration of the algorithm}

To gain more insights into the mechanics
of Algorithm \ref{alg:CRP}, consider the example portrayed in Figure \ref{fig:iteration}. The coordinator begins with $y^A$ as the initial value of $y$. It communicates with SO$_i$ to obtain the critical region induced by $J_i^*$ containing $y^A$, and the affine description of $J_i^*$ over that critical region. Using the relation in \eqref{eq:alphaBeta}, it then computes the critical region $\Pcal^A$ induced by the aggregate cost and the affine description of that cost $\[\alpha^A\]^\T z + \beta^A$ over that region. For convenience, we use
$$ \Pcal^A:= \Pcal^{y^A}, \quad \alpha^{A} := \alpha^{y^A}, \quad  \beta^{A} := \beta^{y^A},$$
and extend the corresponding notation for $y^B, \ldots, y^E$.

\psset{unit=0.38in}
\begin{figure}[h!]
\centering
\vspace{-0.5cm}
\begin{pspicture}[showgrid=false](-2, -1)(8, 5)
        \pstGeonode[PointName=none, PointSymbol=none]
            (0,0){TA1}
            (1,4){TA2}
            (5,4){TA3}
            (6,0){TA4}
            (2,4){TB1}
            (4,0){TB2}
            (0.25,1){TB3}
            (5.25,3){TB4}

    \pstInterLL[PointName=none,PointSymbol=none]{TB1}{TB2}{TB3}{TB4}{TC}

    \psline (TA1)(TA2)(TA3)(TA4)(TA1)

    \pspolygon[fillstyle=solid, fillcolor=darkYellowGreen](TA2)(TB1)(TC)(TB3)
    \pspolygon[fillstyle=solid, fillcolor=SeaGreen](TB1)(TA3)(TB4)(TC)
    \pspolygon[fillstyle=solid, fillcolor=darkOliveGreen](TC)(TB4)(TA4)(TB2)
    \pspolygon[fillstyle=solid, fillcolor=Green](TB3)(TC)(TB2)(TA1)

    \pstGeonode[PointName={\color{White}y^A}, PointSymbol = none, linecolor=White, fillcolor=White, PosAngle = 70, PointNameSep=1.3em](4.5, 1.5){A}
    \pstCircleOA[linecolor=White, fillcolor=White, fillstyle=solid, Radius=1.8]{A}{}

    \pstGeonode[PointName={\color{White}y^B}, PointSymbol = none, linecolor=myRed, fillcolor=myRed, PosAngle = 60, PointNameSep=1.5em](TC){B}
    \pstCircleOA[linecolor=Yellow, fillcolor=myRed, fillstyle=solid, Radius=2.0]{B}{}

    \pstGeonode[PointName={\color{White}y^C}, PointSymbol = none, linecolor=White, fillcolor=White, PosAngle = 90, PointNameSep=1.0em](1.5, 2.5){C}
    \pstCircleOA[linecolor=White, fillcolor=White, fillstyle=solid, Radius=1.8]{C}{}

    \pstGeonode[PointName={\color{White}y^D}, PointSymbol = none, linecolor=White, fillcolor=White, PosAngle = 0, PointNameSep=1.2em](2.5, 1.0){D}
    \pstCircleOA[linecolor=White, fillcolor=White, fillstyle=solid, Radius=1.8]{D}{}

    \pstGeonode[PointName= y^E, PointSymbol =none, linecolor=myRed, fillcolor=myRed, PosAngle = 200, PointNameSep=1.2em](TA1){E}
    \pstCircleOA[linecolor=Yellow, fillcolor=myRed, fillstyle=solid, Radius=2.0]{E}{}

    \pcline[linewidth=1pt, arrowscale=1.5, linecolor=newColorLine, nodesep=3pt]{->}(A)(B)
    \pcarc[linewidth=1pt, arcangle=15, arrowscale=1.5, linecolor=newColorLine, nodesep=3pt]{->}(B)(C)
    \pcarc[linewidth=1pt, arcangle=15, arrowscale=1.5, linecolor=newColorLine, nodesep=3pt]{->}(C)(B)
    \pcline[linewidth=1pt, arrowscale=1.5, linecolor=newColorLine, nodesep=3pt]{->}(B)(D)
    \pcline[linewidth=1pt, arrowscale=1.5, linecolor=newColorLine, nodesep=3pt]{->}(D)(E)

    \rput[bl](6.5, 0.5){\rnode{PA}{\Large $\mathcal{P}^A$}}
    \rput[bl](-1.0, 3.2){\rnode{PC}{\Large $\mathcal{P}^C$}}
    \rput[bl](-1.0, 0.6){\rnode{PD}{\Large $\mathcal{P}^D$}}

    \pcarc[linewidth=2pt, arcangle = 15, arrowscale=1.5, linecolor=Black, nodesep=15pt]{->}(PA)(4.6, 0.7)
    \pcarc[linewidth=2pt, arcangle = -15, arrowscale=1.5, linecolor=Black, nodesep=10pt]{->}(PC)(1.8,3.4)
    \pcarc[linewidth=2pt, arcangle = -15, arrowscale=1.5, linecolor=Black, nodesep=10pt]{->}(PD)(1.2,0.8)

\end{pspicture}
\caption{An example to illustrate the iterative process of Algorithm \ref{alg:CRP}.} \label{fig:iteration}
\end{figure}
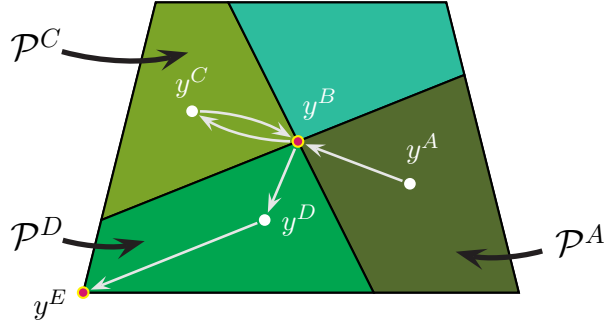

The coordinator solves a linear program to minimize the affine aggregate cost $\[\alpha^A\]^\T z + \beta^A$ over $z \in \Pcal^A$, and obtains the lexicographically smallest optimizer $\yp$.
Such an optimizer $\yp$ is always a vertex of $\Pcal^A$. Identify $y^B$ as that vertex in Figure \ref{fig:iteration}. The optimal cost at $y^B$ is indeed lower than the initial value of $J^*=\infty$, and hence, the coordinator sets $y^* \gets y^B$. It also updates $J^*$ to the aggregate cost at $y^B$, and the partial sub-differential set to $\Dcal \gets \{\alpha^A\}$.

Next, the coordinator solves the least square problem described in \eqref{eq:ifopt} to compute $v^*$. In so doing, it utilizes $\Dcal = \{\alpha^A\}$, and $\v{K}^y = 0$ that describes the normal cone to $\Ycal$ at $y_B$.\footnote{The normal cone to $\Ycal$ at $y^B$ is $\{ 0 \}$ because $y^B$ lies in the interior of $\Ycal$.} Suppose $v^* \neq 0$. The coordinator updates the value of $y$ to $y^C$, obtained by moving a `small' step of length $\ve$ from $y^B$ along $-v^*$. Recall that $y^C \notin \Pcal^A$. The coordinator again communicates with the SOs to obtain the new critical region $\Pcal^C$ induced by the aggregate cost that contains $y^C$. Again, it obtains the affine description of that cost and optimizes it over $\Pcal^C$ to obtain the new $\yp$. In the figure, we depict the case when $\yp$ coincides with $y^*= y^B$.

Notice that the optimal cost $\Jp$ at $\yp$ is equal to $J^*$, and hence, the coordinator only updates the partial sub-differential set $\Dcal$ to $\{ \alpha^A, \alpha^C \}$. With the updated set of $\Dcal$, the coordinator solves \eqref{eq:ifopt} to obtain $v^*$. In this example, $v^*$ is again non-zero, and hence, the coordinator moves along a step of length $\ve$ along $-v^*$ from $y^B$ to land at $y^D$. Again, $y^D \notin \{\Pcal^A, \Pcal^C\}$. The coordinator repeats the same steps to optimize the aggregate cost over $\Pcal^D$ to obtain $y^E$ as the new $\yp$. Two cases can now arise, that we describe separately.
\begin{description}[font=$\bullet$\scshape\space\normalfont, leftmargin=0.2cm]
\item If the optimal cost $\Jp$ at $\yp = y^E$ does not improve upon the cost $J^*$ at $y^B$, the coordinator ignores $y^E$ and updates the set $\Dcal$ to $\{ \alpha^A, \alpha^C, \alpha^D \}$. It computes $v^*$ with the updated $\Dcal$. Again, if $v^* \neq 0$, it traverses along $-v^*$ to venture into a yet-unexplored critical region. The process continues till we get $y^*=y^B$ as an optimizer (if $v^* = 0$ at a future iterate), or we encounter the case we describe next.

\item If $\Jp < J^*$, then the coordinator sets $y^E$ as the new $y^*$. It retraces the same steps with this new $y^*$. In this example, since $y^E$ is a vertex of $\Ycal$, one can show that  \eqref{eq:ifopt} will yield $v^* = 0$, and hence, $y^* = y^E$ will optimize the aggregate cost over $\Ycal$.
\end{description}



\newcommand{\xiL}{\xi^{\sf L}}
\newcommand{\xiU}{\xi^{\sf U}}

\section{The robust counterpart}
\label{sec:robust}

The deterministic tie-line scheduling problem was formulated in the last section on the premise that available  renewable supply and limits on power demands within each area are known at the time when tie-line schedules are decided. We now alter that assumption and allow these parameters to be uncertain. In particular, we let $\xi_i = \(\overline{W}_i, \underline{D}_i, \overline{D}_i \)$ take values in a \emph{box}, described by
\begin{equation}
\Xi_i := \{ \xi_i \in \Rset^{3 n_i} \ \vert \ \xiL_i \leq \xi_i \leq \xiU_i \}\label{eq:Xi}
\end{equation}
for $i=1,2$. The robust counterpart of the tie-line scheduling problem is then described by
\begin{align}
\underset{y \in \Ycal}{\text{minimize}} \( \underset{\xi_1 \in \Xi_1}{\text{max}} \ J_1^*\(y, \xi_1\) + \underset{\xi_2 \in \Xi_2}{\text{max}} \ J_2^*\(y, \xi_2\) \).
\label{eq:robProb}
\end{align}
We now develop an algorithm that solves \eqref{eq:robProb} in a distributed fashion.
Problem \eqref{eq:robProb} has a \emph{minimax} structure. Therefore, we employ a strategy in Algorithm \ref{alg:robust} to alternately minimize the objective function over $\Ycal$ and maximize it over $\Xi_1 \times \Xi_2$. Thanks to the following lemma, the maximization over $\Xi_1 \times \Xi_2$ can be reformulated into a mixed-integer linear program.
\begin{lemma}
\label{lemma:MILP}
Fix $y \in \Ycal$. Then, there exists $\Msf > 0$ for which maximizing $J_i^*(y, \xi_i)$ over $\xi_i \in \Xi_i$ is equivalent to the following mixed-integer linear program:
\begin{equation} \label{eq:MILP}
\hspace{-0.3cm}
\begin{aligned}
& \underset{w_i, \rho, \lambda}{\text{maximize}}
	 && c_i^0 + [c_i^\xi]^\T \ \xiL_i \!+\! ( \v{A}_i^\xi \xiL_i \! + \! \v{A}_i^y y \!-\! b_i)^\T \ \lambda \!+\!  \bone^\T \rho, \\
 & \text{subject to}   \quad
	&& c_i^x + [\v{A}_i^x]^\T \lambda = 0, \\ 
	&&& \rho \leq {\Msf} w_i, \\
	&&& \rho \leq {\Msf} (\bone - w_i) + \v{\Delta}^\xi_{i}(c_i^\xi + [\v{A}_i^\xi ]^\T \lambda), \\
	&&& w_i \in \{ 0, 1 \}^{n_i}, \rho \in \Rset^{n_i}, \lambda \in \Rset_+^{m_i}.
 \end{aligned}\!\!
\end{equation}
\end{lemma}
We use the notation $\v{\Delta}^\xi_{i}$ to denote a diagonal matrix with $\xiU_i-\xiL_i$ as the diagonal.
The lemma builds on the fact that $J_i^*(y, \xi_i)$ is convex in $\xi_i$, and hence, reaches its maximum at a vertex of $\Xi_i$. The convexity is again a consequence of \cite[Theorem 7.5]{borrelli2003constrained}. Our proof in Appendix \ref{sec:lemma2} leverages duality theory of linear programming and the so-called \emph{big-M} method adopted in \cite[Chapter 2.11]{conforti2014integer} to reformulate the maximization of $J_i^*(y, \cdot)$ over the vertices of $\Xi_i$ into a mixed-integer linear program. An optimal $\xi_i^{\text{opt}}$ can be recovered from $w_i^*$ that is optimal in \eqref{eq:MILP} using
$$ \xi_i^{\text{opt}} := \xiL_i + \v{\Delta}^\xi_{i} w_i^*.$$

Next, we present our algorithm for solving the robust counterpart. In the algorithm, the SOs exclusively maintain and update certain variables; we distinguish these from the ones the coordinator maintains.
\begin{algorithm}
	\caption{Solving the robust counterpart.}
	\label{alg:robust}
	\begin{algorithmic}[1]
		\Initialize{
		SO$_1$: $\Vcal_1 \gets \{ \text{a vertex of }\Xi_1 \}$, \\
		SO$_2$: $\Vcal_2 \gets \{ \text{a vertex of }\Xi_2 \}$.}
		\Do
	        		\State Coordinator uses Algorithm \ref{alg:CRP} to solve
			$$\underset{y \in \Ycal}{\text{minimize}} \ \(\max_{\xi_1 \in \Vcal_1} J_1^*(y, \xi_1) + \max_{\xi_2 \in \Vcal_2} J_2^*(y, \xi_2)\). $$\label{step:CRP}
			\State $y^* \gets$ optimizer in step \ref{step:CRP}.
			\State $J^* \gets$ optimal cost in step \ref{step:CRP}.
			\State For $i=1,2$, SO$_i$ performs:
			\State \qquad Maximize $J_i^*(y^*, \cdot)$ over $\Xi_i$ using \eqref{eq:MILP}.\label{step:MILP}
			\State \qquad $\xi_i^{\text{opt}} \gets$ optimizer in step \ref{step:MILP}.
			\State \qquad $\Jp_i \gets$ optimal cost in step \ref{step:MILP}.

			\State \qquad $\Vcal_i \gets \Vcal_i \cup \{ \xi_i^\text{opt} \}$.
			\State \qquad \Return $\Jp_i$ to the coordinator.
		\doWhile{$\Jp_1 + \Jp_2 > J^*$.}
	\end{algorithmic}
\end{algorithm}

We summarize the main property of the above algorithm in the following theorem, whose proof is given in Appendix \ref{sec:theorem2}\footnote{The proof is similar to \cite[Preposition 2]{zeng2013solving}; we include it for completeness.}.
\begin{theorem}
\label{thm:finiteTimeRob}
Algorithm \ref{alg:robust} terminates after finitely many steps, and $y^*$ at termination optimally solves \eqref{eq:robProb}.
\end{theorem}

Our algorithm to solve the robust counterpart makes use of Algorithm \ref{alg:CRP} in step \ref{step:CRP}. The coordinator performs this step with necessary communication with the SOs. However, it remains agnostic to the uncertainty sets $\Xi_1$ and  $\Xi_2$ throughout. Therefore, our algorithm is such that the SOs in general will not be required to reveal their cost structures, network constraints, nor their uncertainty sets to the coordinator to optimally solve the robust tie-line scheduling problem. Further, Theorems \ref{thm:finiteTime} and \ref{thm:finiteTimeRob} together guarantee that the coordinator can arrive at the required schedule by communicating with the SOs only finitely many times. These define some of the advantages of the proposed methodology. In the following, we discuss some limitations of our method.

The number of affine segments in the piecewise affine description of
$\max_{\xi_i \in \Vcal_i} J_i^*(y, \xi_i)$
increases with the size of the set $\Vcal_i$. The larger that number, the heavier can be the computational burden on Algorithm \ref{alg:CRP} in step \ref{step:CRP}. To partially circumvent this problem, we initialize the sets $\Vcal_i$ with that vertex of $\Xi_i$ that encodes the least available renewable supply and the highest nominal demand. Such a choice captures the intuition that dispatch cost is likely the highest with the least free renewable supply and the highest demand. Our empirical results in the next section corroborate that intuition.

We make use of mixed-integer linear programs in step \ref{step:MILP} of the algorithm. This optimization class encompasses well-known NP-hard problems. Solvers in practice, however, often demonstrate good empirical performance. Popular techniques for mixed-integer linear programming include branch-and-bound, cutting-plane methods, etc. See \cite{conforti2014integer} for a survey. Providing polynomial-time convergence guarantees for \eqref{eq:MILP} remains challenging, but our empirical results in the next section appear encouraging.


\section{Numerical Experiments}
\label{sec:simulation}

We report here the results of our implementation of Algorithm \ref{alg:robust} on several power system examples. All optimization problems were solved in IBM ILOG CPLEX Optimization Studio V12.5.0 \cite{CPLEX} on a PC with 2.0GHz Intel(R) Core(TM) i7-4510U microprocessor and 8GB RAM.

\begin{figure}[h!]
    \centering
    \begin{subfigure}[t]{0.45\textwidth}
        \centering
        \includegraphics[width=\textwidth]{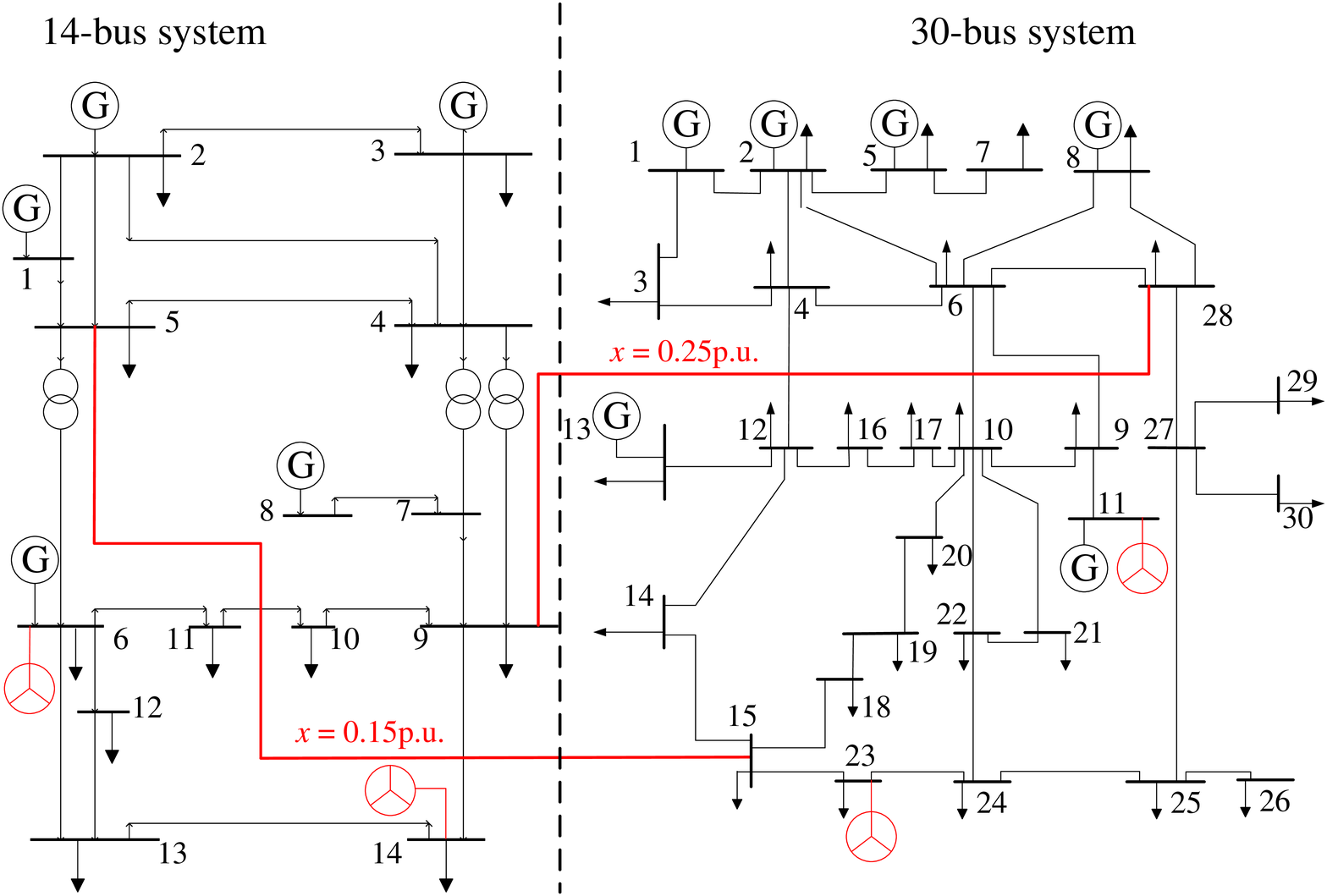}
        \caption{The power system model}
        \label{fig:twoareaConfig}
    \end{subfigure}
    \begin{subfigure}[t]{0.46\textwidth}
        \centering
        \includegraphics[width=\textwidth]{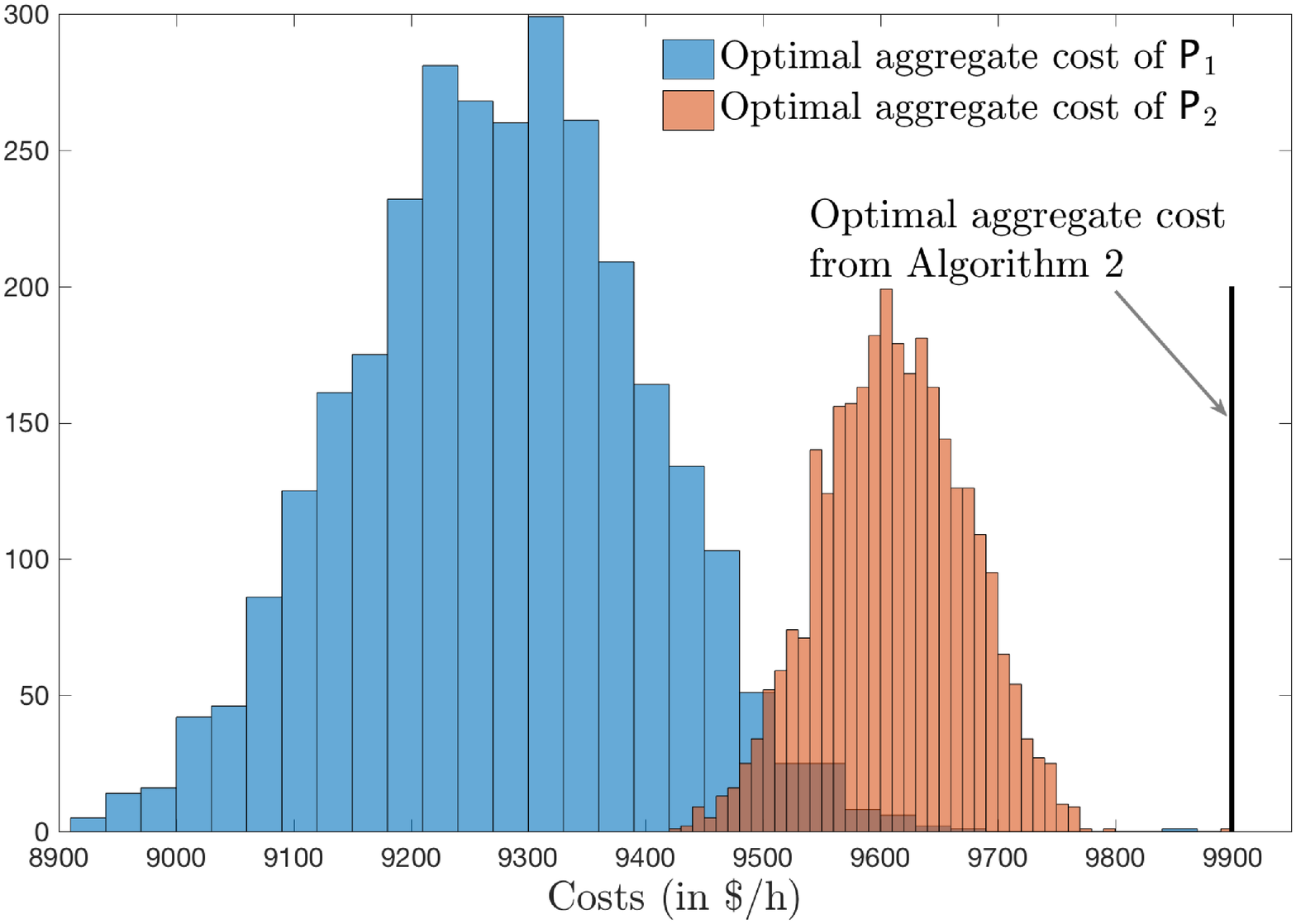}
        \caption{Histogram of optimal aggregate costs}
        \label{fig:twoareaCost}
    \end{subfigure}
    \caption{The two-area 44-bus system is portrayed on the left. It shows where the wind generators are added and the parameters for the tie-lines used in our experiments.
    The figure to the right plots the optimal aggregate costs from $\Psf_1$, $\Psf_2$ over 3000 samples of uncertain variables, and that of Algorithm \ref{alg:robust} on this system.}
\end{figure}
\subsection{On a two-area 44-bus power system}
\label{sec:2area}

Consider the two-area power system shown in Figure \ref{fig:twoareaConfig}, obtained by connecting the IEEE 14- and 30-bus test systems \cite{IEEESys}. The networks were augmented with wind generators at various buses. Transmission capacities of all lines were set to 100MW. The available capacity of each wind generator was varied between 15MW and 25MW. The lower limits on all power demands were set to zero, while the upper limits were varied between 98\% and 102\% of their nominal values. Our setup had 36 uncertain variables -- 32 power demands and 4 available wind generation. Bus 5 in area 1 was the slack bus.
From the data in Matpower \cite{MATPOWER}, we chose the linear coefficient in the nominal quadratic cost structure for each conventional generator to define $P^g_i$ in \eqref{eq:defP}. Further, we neglected wind spillage costs by letting $P^w_i = 0$, and defined $P^d_i$ by assuming a constant marginal cost of \$100/MWh for not meeting the highest demands.

\begin{table}[H]
\centering
\begin{tabular}{@{\hskip 0.05in}c@{\hskip 0.05in}l@{\hskip 0.05in}c@{\hskip 0.05in}c@{\hskip 0.05in}}
	\hlinewd{1pt}
  	Iteration & Step in Algorithm \ref{alg:robust}  &  \thead{Aggregate cost \\ (in \$/h)} &  \thead{Run-time \\ (in ms)}   \\
  	\hline
          1 & Step \ref{step:CRP} to compute $y^*$    & 9897.7 & 113.6\\
          1 & Step \ref{step:MILP} to compute $\xi^{\text{opt}}$ & 9910.3  & 99.6\\
          2 & Step \ref{step:CRP} to compute $y^*$    & 9899.3  & 93.4\\
          2 & Step \ref{step:MILP} to compute $\xi^{\text{opt}}$ & 9899.3 & 121.5\\
	\hlinewd{1pt}
\end{tabular}
\caption{\small Evolution of aggregate cost of Algorithm \ref{alg:robust} for the two-area power system in Figure \ref{fig:twoareaConfig}.}
\label{table:twoarearesult}
\end{table}


To run Algorithm \ref{alg:robust}, we initialized $\Vcal_i$ with the scenario that describes the highest power demands and the least available wind generation across all buses. To invoke Algorithm \ref{alg:CRP} in step \ref{step:CRP}, we initialized $y$ with a vector of all zeros. When the algorithm encountered the same step in future iterations, it was initialized with the optimal $y^*$ from the last iteration to provide a \emph{warm start}.
Algorithm \ref{alg:robust} converged in two iterations, \textit{i.e.}, it ended when the cardinality of $\Vcal_1$ and $\Vcal_2$ were both two. The trajectory of the optimal cost and the run-times for each step are given in Table \ref{table:twoarearesult}. In the first iteration, Algorithm \ref{alg:CRP} in step \ref{step:CRP} with $\ve = 10^{-5}$ converged in four iterations\footnote{The termination condition $ v^* = 0 $ is replaced by checking that the Euclidean norm of a suitably normalized $v^*$ is less than a threshold.} of its own and explored five critical regions induced by the aggregate cost. A naive search over $\Ycal$ yielded that the aggregate cost induced at least 126 critical regions. Our simulation indicates that Algorithm \ref{alg:CRP} only explores a `small' subset of all critical regions.

Step \ref{step:MILP} of Algorithm \ref{alg:robust} was then solved to obtain $\xi_i^{\text{opt}}$. As Table \ref{table:twoarearesult} suggests, the aggregate cost $J_1^{\text{opt}} + J_2^{\text{opt}}$ exceeded $J^*$ obtained earlier in step \ref{step:CRP}. Thus, the scenario of demand and supply captured in our initial sets $\Vcal_1$ and $\Vcal_2$ was \emph{not} the one with maximum aggregate dispatch costs. To accomplish this step, two separate mixed-integer linear programs were solved --  one with 13 binary variables (in area 1) and the other with 23 binary variables (in area 2). CPLEX returned the global optimal solutions in 15ms and 77ms, respectively.
In the next iteration, step \ref{step:CRP} was performed with $\xi_i^{\text{opt}}$ added to $\Vcal_i$, where Algorithm \ref{alg:CRP} converged in five iterations, exploring only four critical regions. Finally, step \ref{step:MILP} yielded $\Jp_1 + \Jp_2 = J^*$, implying that the obtained $y^*$ defines an optimal robust tie-line schedule.

To further understand the efficacy of our solution technique, we uniformly sampled the set $\Xi_1 \times \Xi_2$ 3000 times. With each sample $\(\xi_1, \xi_2\)$, we solved two optimization problems -- $\Psf_1$ and $\Psf_2$. Precisely, $\Psf_1$ is a deterministic tie-line scheduling problem solved with Algorithm \ref{alg:CRP}, and $\Psf_2$ is the optimal power flow problem in each area with the optimal $y^*$ obtained from Algorithm \ref{alg:robust} for the robust counterpart. The histograms of the optimal aggregate costs from $\Psf_1$ and $\Psf_2$ are plotted in Figure \ref{fig:twoareaCost}. The same figure also depicts the optimal cost of the robust tie-line scheduling problem, which naturally equals the maximum among the costs from $\Psf_2$. And for each sample, the gap between the optimal costs of $\Psf_1$ and $\Psf_2$ captures the cost due to lack of foresight. Figure \ref{fig:twoareaCost} reveals that such costs can be significant. The median run-time of $\Psf_1$ was 48.5ms over all samples. The run-time for the robust problem was 458.2ms -- roughly 10 times that median.

\subsection{On a three-area 187-bus system test}
\label{sec:3area}
\begin{figure}[h!]
    \centering
    \begin{subfigure}[t]{0.48\textwidth}
        \centering
        \includegraphics[width=\textwidth]{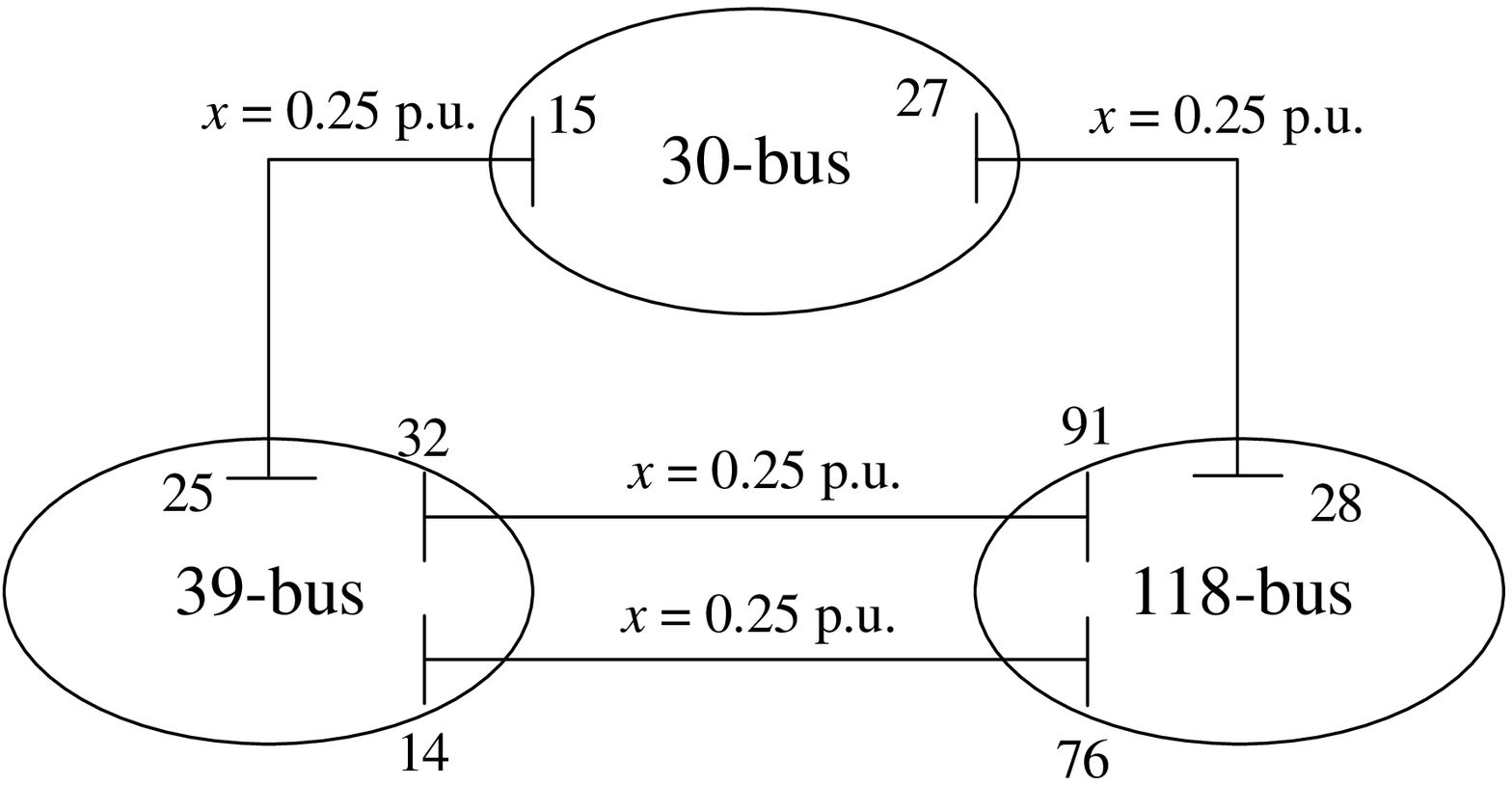}
        \caption{A three-area 187-bus power system.}
        \label{fig:threeareaConfig}
    \end{subfigure}
    \quad
    \begin{subfigure}[t]{0.48\textwidth}
        \centering
        \includegraphics[width=\textwidth]{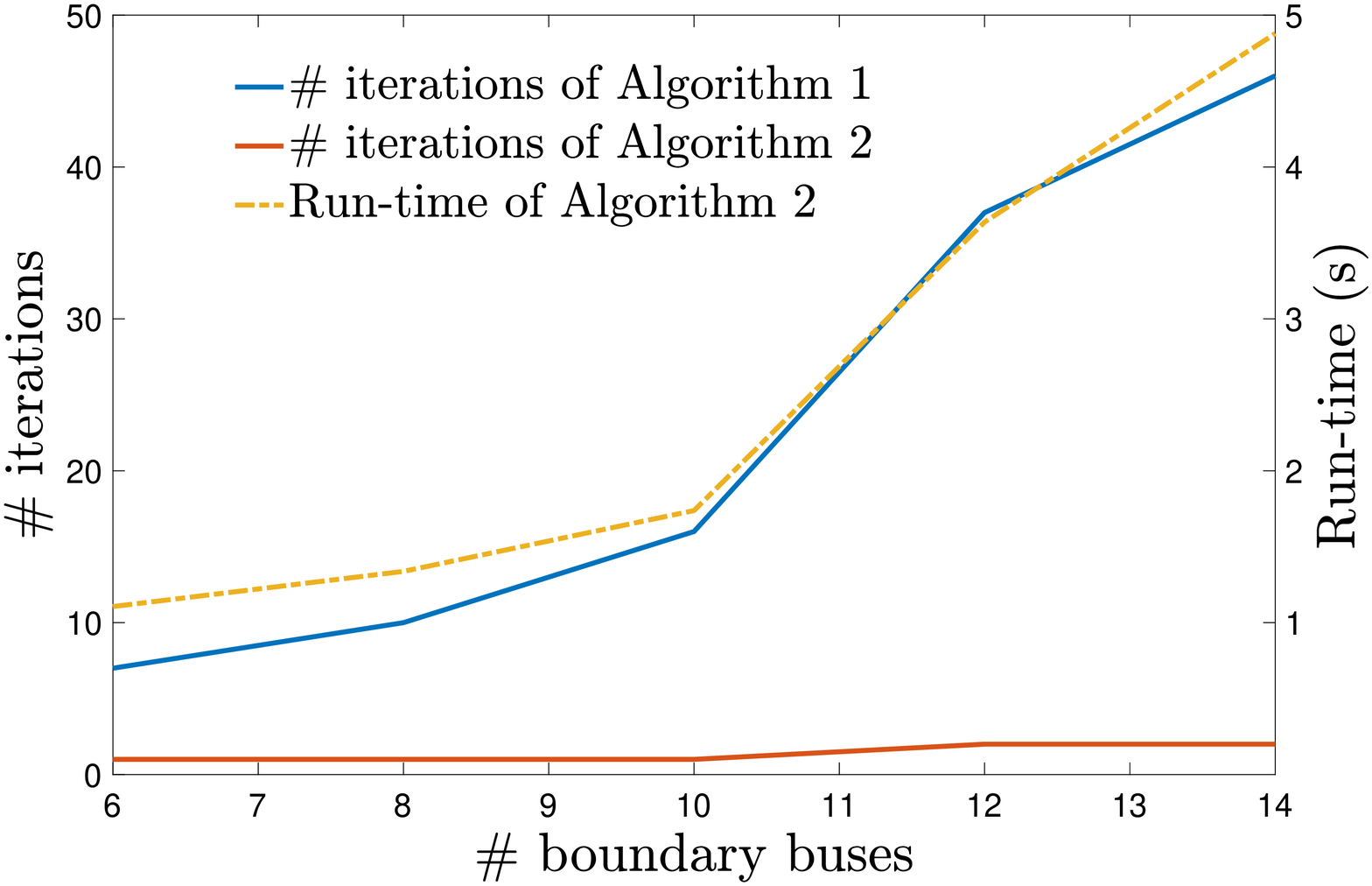}
        \caption{Performance of algorithms with $\#$ of tie-lines.}
        \label{fig:Perturbation}
    \end{subfigure}
    \caption{The power system model and how our algorithms perform with variation in number of tie-lines in a three-area 187-bus power system.}
\end{figure}

For this case study, we interconnected the IEEE 30-, 39-, and 118-bus test systems as shown in Figure \ref{fig:threeareaConfig}. All transmission capacities were set to 100MW. Five wind generators were added to the 118-bus system (at buses 17, 38, 66, 88, and 111), three in the 39-bus system (at buses 3, 19, and 38), and two in the 30-bus system (at buses 11, and 23). Again, we adopted the same possible set of available wind power generations and power demands, as well as the cost structures as in Section \ref{sec:2area}. In total, our robust tie-line scheduling problem modeled 151 uncertain variables. For this multi-area power system, Algorithm \ref{alg:robust} converged in the first iteration. The mixed integer programs in step \ref{step:MILP} yielded the global optimal solution for each area, taking 62ms, 109ms, and 281ms, respectively.
We again sampled the set $\Xi_1 \times \Xi_2 \times \Xi_3$ 3000 times, and solved $\Psf_1$.
The run-time of Algorithm \ref{alg:robust} was 825.3ms, that is roughly 1.8 times the median run-time of $\Psf_1$, given by 450.8ms.

We studied how our algorithm scales with the number of boundary buses by adding more tie-lines to the same system. The aggregate iteration count of Algorithm \ref{alg:CRP} is expected to grow with the number of induced critical regions, that in turn should grow with the boundary bus count. On the other hand,
the iteration count of Algorithm \ref{alg:robust} largely depends on the initial choice of the scenario encoded in the sets $\Vcal_1, \Vcal_2, \Vcal_3$, and thus, varies to a lesser extent on the same count. Figure \ref{fig:Perturbation} validates these intuitions.

\subsection{Summary of results from other case-studies}
\label{sec:othertest}

We compared Algorithm \ref{alg:CRP} with a dual decomposition based approach proposed in \cite{Bakirtzis&Biskas:03TPS}. That algorithm converges asymptotically, while our method converges in finitely many iterations. Table \ref{table:comparison} summarizes the comparison.\footnote{We say the method in \cite{Bakirtzis&Biskas:03TPS} converges when the power flow over each tie-line as calculated by the areas at its end mismatches by $<$ 0.01 p.u..} Compared to that in \cite{Bakirtzis&Biskas:03TPS}, our algorithm clocked lesser number of iterations and lower run-times in our experiments.

\begin{table}[H]
		\centering
		\begin{tabular}{lcccc}
                    	\hlinewd{1pt}
                      	Items & \thead{Two-area \\44-bus system}  &  \thead{Three-area\\ 187-bus system}  \\
                      	\hline
                      		{\# iterations in Algorithm \ref{alg:CRP}} & 8  & 9 \\
                      		{\# iterations of \cite{Bakirtzis&Biskas:03TPS}} & 23  & 78 \\
                      		{Run-time of Algorithm \ref{alg:CRP} (ms)} & 458.2  & 825.3 \\
                      		{Run-time of \cite{Bakirtzis&Biskas:03TPS} (ms)}  & 779.8  & 1227.5 \\
                  	\hlinewd{1pt}
        		\end{tabular}
		\caption{Comparison with the method in \cite{Bakirtzis&Biskas:03TPS}.}%
		\label{table:comparison}
\end{table}

Apart from the two systems considered so far, we ran Algorithm \ref{alg:robust} on a collection of other multi-area power systems, details of which can be found in Appendix \ref{app:otherSystems}.
The results are summarized in  Table \ref{table:moreresults}. Our experiments reveal that Algorithm \ref{alg:robust} often converges within 1 -- 4 iterations.
The run-time of Algorithm \ref{alg:robust} grows significantly with the number of uncertain parameters. The 418-bus and the 536-bus systems with 422 and 546 uncertain variables, respectively, corroborate that conclusion. Such growth in run-time is expected because the complexity of \eqref{eq:MILP} grows with the number of binary decision variables that equals the number of uncertain parameters. Run-time of a joint multi-area optimal power flow problem with a sample scenario in the last column provides a reference to compare run-times for the robust one.

\begin{table}[H]
		\centering
        		\begin{tabular}{ccccccc}
              	\hlinewd{1pt}
                      \thead{\# \\ areas} & \thead{\# \\ buses}  & \thead{\# \\ uncertain\\ variables} & \thead{\# \\ boundary\\ buses} &
                      \thead{\# iter. in \\ Algorithm \ref{alg:robust}} &
                      \thead{Run-time of \\ Algorithm \ref{alg:robust}\\ (ms)} &
                      \thead{Run-time of \\ joint problem \\ (ms)} \\
                      \hline
                      2   & 87        	& 91          &4& 		1                  	& 719.6   	&  310.0 \\
                      2   & 175       	& 179         &4& 		 1                 	& 871.1 	&  340.5 \\
                      2   & 236       	& 240         &4&		 1                  	& 1732.6 	&  391.5 \\
                      2   & 418       	& 42          &10&  		 1                  	& 1020.7  	&  455.7 \\
                      2   & 418       	& 422         &10& 		 4                  	& 6124.5  &  461.4 \\
                      3   & 354       	& 360         &12& 		 3                  	& 4127.4 	&  655.8   \\
                      3   & 536       	& 54          &12&  	 1                  	& 2557.6 	&  699.7  \\
                      3   & 536       	& 546         &12& 		 3                  	& 18359.8 &  701.2 \\
                      \hlinewd{1pt}
        		\end{tabular}
		\caption{Performance of Algorithm \ref{alg:robust} on various multi-area power system examples provided in Appendix \ref{app:otherSystems}.}%
		\label{table:moreresults}
\end{table}


\section{Conclusion}
\label{sec:conclusion}

This work presented an algorithmic framework to solve a tie-line scheduling problem in multi-area power systems. Our method requires a coordinator to communicate with the system operators in each area to arrive at an optimal tie-line schedule. In the deterministic setting, where the demand and supply conditions are assumed known during the scheduling process, our method (Algorithm \ref{alg:CRP}) was proven to converge in finitely many steps. In the case with uncertainty, we proposed a method (Algorithm \ref{alg:robust}) to solve the robust variant of the tie-line scheduling problem. Again, our method was shown to converge in finitely many steps. Our proposed algorithms do not require the system operator to reveal the dispatch cost structure, network parameters or even the support set of uncertain demand and supply within each area to the coordinator. We empirically demonstrated the efficacy of our algorithms on various multi-area power system examples.


\bibliographystyle{plain}
\bibliography{ref}

\appendix

\section{How SO$_i$ can compute $\Pcal_i^y$, $\alpha_i^y$, $\beta_i^y$}
\label{sec:CRAffine}

With $\xi_i\in\Xi_i$ and $y \in \Ycal$ fixed, consider the optimization problem described in \eqref{eq:areaProb}. Suppose the optimal solution $x_i^*(y, \xi_i)$ is unique. We suppress the dependency on $(y, \xi_i)$ for notational convenience. Distinguish between the constraints that are \emph{active} (met with an equality) versus that are \emph{inactive} at optimality with the subscript $\Acal$ and $\Ical$, respectively, as follows.
\begin{align*}
[\v{A}_i^x]_{\Acal} \ x_i^* + [\v{A}_{i}^\xi]_{\Acal} \ \xi_i + \[\v{A}_{i}^y\]_{\Acal}  y &=  [b_i]_{\Acal},\\
[\v{A}_i^x]_{\Ical} \ x_i^* + [\v{A}_{i}^\xi]_{\Ical} \ \xi_i + \[\v{A}_{i}^y\]_{\Ical} y &<  [b_i]_{\Ical}.
\end{align*}
The set of active versus inactive constraints remains the same over the critical region $\Pcal_i^y$. Assuming $[\v{A}_i^x]_{\Acal}$ is a square and invertible matrix, the optimal solution $x_i^*$ is unique for each $z \in \Pcal^y_i$, given by
\begin{align*}
x_i^* = [\v{A}_i^x]_{\Acal}^{-1} \( [b_i]_{\Acal} - [\v{A}_{i}^\xi]_{\Acal} \ \xi_i - \[\v{A}_{i}^y\]_{\Acal}  z \).
\label{eq:local_partition1}
\end{align*}
The inequalities for the inactive constraints, together with the above relation defines the critical region $ \Pcal_i^y := \{z \in \Ycal \! : \! \v{D} z \leq d\}$, where
\begin{align*}
\v{D}
&= - [\v{A}_i^x]_{\Ical} [\v{A}_i^x]_{\Acal}^{-1}[\v{A}_{i}^y]_{\Acal}+
[\v{A}_{i}^y]_{\Ical},\\
d
&= [b_i]_{\Ical} - [\v{A}_{i}^\xi]_{\Ical} \ \xi_i  - [\v{A}_i^x]_{\Ical} [\v{A}_i^x]_{\Acal}^{-1} \( [b_i]_{\Acal} -  [\v{A}_{i}^\xi]_{\Acal} \ \xi_i \).
\end{align*}
Finally, $J_i^*(y, \xi_i) = c_i(x_i^*, \xi_i)$ yields
\begin{align*}
\alpha^y_i &= -[c_i^x]^\T \ [\v{A}_i^x]_{\Acal}^{-1}[\v{A}_{i}^y]_{\Acal}, \\
\beta^y_i &= c_i^0 + [c_i^\xi]^\T \ \xi_i +  [c_i^x]^\T \ [\v{A}_i^x]_{\Acal}^{-1} \( [b_i]_{\Acal} -  [\v{A}_{i}^\xi]_{\Acal} \ \xi_i \).
\end{align*}
The above expressions are derived under the premise that $[\v{A}_i^x]_{\Acal}$ is invertible. We refer the reader to \cite[Sections 7.2.2, 7.2.4]{borrelli2003constrained} for the procedure in the general case.


\section{Proof of Theorem \ref{thm:finiteTime}}\label{sec:finiteTime}

After each iteration of Algorithm \ref{alg:CRP}, $y^*$ is a vertex of a critical region induced by the aggregate optimal cost. Also, $\Dcal$ is such that $\conv(\Dcal) \subseteq \partial J^*\(y^*, \xi_1, \xi_2\)$. Therefore, if the algorithm terminates with $v^* =0$, then
$$ 0 \in \conv(\Dcal) + N_\Ycal\(y^*\) \subseteq \partial J^*\(y^*, \xi_1, \xi_2\) + N_\Ycal\(y^*\).$$
That is, $y^*$ optimally solves \eqref{eq:detProb.2}. Next, we argue that the algorithm terminates in finitely many iterations.

Consider the sequence of $y^*$'s and $J^*$'s produced by the algorithm. Notice that $J^*$ is a piecewise constant but non-increasing sequence. Further, a change in $y^*$ always accompanies a strict decrease in $J^*$. Therefore, if $y^*$ changes in an iteration from a certain point, that same point can never become $y^*$ again. Since there are finitely many critical regions with finitely many vertices, it only remains to show that $y^*$ cannot remain constant over infinitely many iterations. Towards that goal, notice that $y^*$ can only belong to a finite number of critical regions. In the rest of the proof, we argue that the variable $y$ computed in step \ref{step:yUpdate} always belongs to a different such critical region containing $y^*$, unless the algorithm terminates.

At an arbitrary iteration, assume that $y$ has taken values in critical regions $\Pcal^1, \ldots, \Pcal^{\ell_D}$ that contain $y^*$. For convenience, let the optimal aggregate cost be given by $\[\alpha^j\]^\T z + \beta^j$ for $z \in \Pcal^j$ for each $j=1,\ldots, \ell_D$. Thus, $\Dcal := \{ \alpha^1, \ldots, \alpha^{\ell_D}\}$.
Then, the new value of $y$ is computed as $y^* - \ve v^*$, with $v^*$ as defined in \eqref{eq:ifopt}. If $v^* = 0$, then the algorithm terminates, proving our claim. Otherwise, assume that $y^* - \ve v^* \in \Pcal^1$, contrary to our hypothesis, implying
\begin{align*}
J^*(y^* - \ve v^*, \xi_1, \xi_2)
= \[\alpha^1\]^\T \(y^* - \ve v^*\) + \beta^1 
= J^* (y^*, \xi_1, \xi_2) - \ve \[ \alpha^1\]^\T v^*.
\end{align*}
Since, $y^*$ optimizes the aggregate cost over $\Pcal^1$, it suffices to show that $\[ \alpha^1\]^\T v^* > 0$ to arrive at a contradiction. For convenience, define the matrix $\pmb{\alpha} := \( \alpha^1, \ldots, \alpha^{\ell_D} \)$. We prove more generally that $\pmb{\alpha}^\T v^* > 0$. Associate Lagrange multipliers $\phi, \psi$ with the equality constraints $v = \pmb{\alpha} \eta + \zeta$, and $\bone^\T \eta = 1$, respectively. Also, associate $\mu_\eta, \mu_\zeta$ with the inequality constraints $\eta \geq 0$ and $\v{K}^y \zeta \geq 0$, respectively. Then, an optimal primal-dual solution pair given by $v^*, \eta^*, \zeta^*$ and $\phi^*, \psi^*, \mu_\eta^*, \mu_\zeta^*$ satisfies the Karush-Kuhn-Tucker (KKT) optimality conditions --  comprised of the constraints in \eqref{eq:ifopt} and the following relations.
\begin{gather*}
{v^*- \phi^*} = 0, \ \
{\pmb{\alpha}^\T \phi^* + \psi^* \bone - \mu_\eta^*} = 0, \ \
{\phi^* - \[\v{K}^y\]^\T \mu_{\zeta}^*} = 0, \\
\[ \mu_{\eta}^*\]^\T \eta = 0, \ \ \[\mu_{\zeta}^*\]^\T \zeta^* = 0, \ \
\mu_{\eta}^* \geq 0, \ \ \mu_\zeta^* \geq 0.
\end{gather*}
Using the KKT conditions, we have
\begin{align*}
\vnorm{v^*}_2^2 + \psi^*
&= \(\phi^*\)^\T\(\pmb{\alpha} \eta^* + \zeta^*\) + \psi^* \\
&= \( \pmb{\alpha}^\T \phi^* + \psi^* \bone - \mu_\eta^* \)^\T \eta^* + \( \phi^* - \[\v{K}^y\]^\T \mu_{\zeta}^* \)^\T \zeta^* \\
&= 0.
\end{align*}
Thus, $\psi^* < 0$. Together with the KKT conditions, that yields
\begin{align*}
\pmb{\alpha}^\T v^* = \pmb{\alpha}^\T \phi^* = -\psi^* \bone + \mu_\eta^* > 0.
\end{align*}


\section{Proof of Lemma \ref{lemma:MILP}}\label{sec:lemma2}
Strong duality of the problem described in \eqref{eq:areaProb} implies that $J_i^*\(y, \xi_i\) $ equals the optimum of the following problem.
\begin{equation*} \label{eq:areaDual}
\hspace{-0.4cm}
\begin{alignedat}{8}
& \quad \underset{\lambda \in \Rset^{m_i}_+}{\text{maximum}}
	 & & c_i^0 + [ c_i^\xi ]^\T \xi_i +  \( \v{A}_{i}^y y +  \v{A}_{i}^{\xi} \xi_i - b_i \)^\T \lambda, \\
 & \quad \text{subject to}   \quad
& \quad & c_i^x + \[ \v{A}_i^x \]^\T \lambda = 0.
 \end{alignedat}
\end{equation*}
Then, maximizing $J_i^*\(y, \xi_i\)$ over the vertices of $\Xi_i$, described by
$ \{ \xiL_i + \v{\Delta}^\xi_i w_i \! : \! w_i \in \{0, 1\}^{n_i} \}$, is equivalent to
\begin{equation*}
\hspace{-0.4cm}
\begin{alignedat}{8}
& \quad {\text{maximize}}
	 & & c_i^0 + [c_i^\xi]^\T \xiL_i \!+\! \( \v{A}_i^y y  \! + \! \v{A}_i^\xi \xiL_i \!-\! b_i\)^\T \lambda \!+\!  \bone^\T \rho, \\
 & \quad \text{subject to}   \quad
& \quad & c_i^x + \[ \v{A}_i^x \]^\T \lambda = 0, \\
&&& \rho = \diag(w_i) \cdot  \v{\Delta}_i^\xi \( c_i^\xi + [ \v{A}_i^\xi ]^\T \lambda  \).
\end{alignedat}
\end{equation*}
over $w_i \in \{0, 1\}^{n_i}$, $\rho\in\Rset^{n_i}$, and $\lambda \in \Rset^{m_i}_+$.
Here, $\diag(w_i)$ denotes the diagonal matrix with $w_i$ as the diagonal. Since we maximize $\bone^\T \rho$, one can replace the second equality constraint in the above problem with the inequality
$$\rho \leq \diag(w_i) \cdot \v{\Delta}_i^\xi \( c_i^\xi + [ \v{A}_i^\xi ]^\T \lambda  \),$$
that is further equivalent to
\begin{align*}
\rho \leq {\Msf} w_i, \ \text{and} \ \  \rho \leq {\Msf} \(\bone - w_i\) + \v{\Delta}_i^\xi \( c_i^\xi + [ \v{A}_i^\xi ]^\T \lambda  \),
\end{align*}
for a large enough ${\Msf} > 0$. That completes the proof.


\section{Proof of Theorem \ref{thm:finiteTimeRob}}\label{sec:theorem2}
Let $J^{\text{rob}}$ denote the optimal aggregate cost of \eqref{eq:robProb}. Then, $J^*$ from step \ref{step:CRP} and $\Jp_1 + \Jp_2$ from step \ref{step:MILP} at any iteration of Algorithm \ref{alg:robust} satisfy
$$ J^* \leq J^{\text{rob}} \leq \Jp_1 + \Jp_2.$$
If Algorithm \ref{alg:robust} terminates, the termination condition implies that the above inequalities are all equalities. In that event, $y^*$ optimally solves \eqref{eq:robProb}.

To argue the finite-time convergence, notice that at least one among $\Vcal_1$ and $\Vcal_2$ increases in cardinality unless the termination condition is satisfied. The rest follows from the fact that $\Xi_1$ and $\Xi_2$ have finitely many vertices.


\section{Power system details for additional simulations}
\label{app:otherSystems}

\begin{figure}[h!]
    \centering
    \begin{subfigure}[t]{0.35\textwidth}
        \centering
        \includegraphics[width=\textwidth]{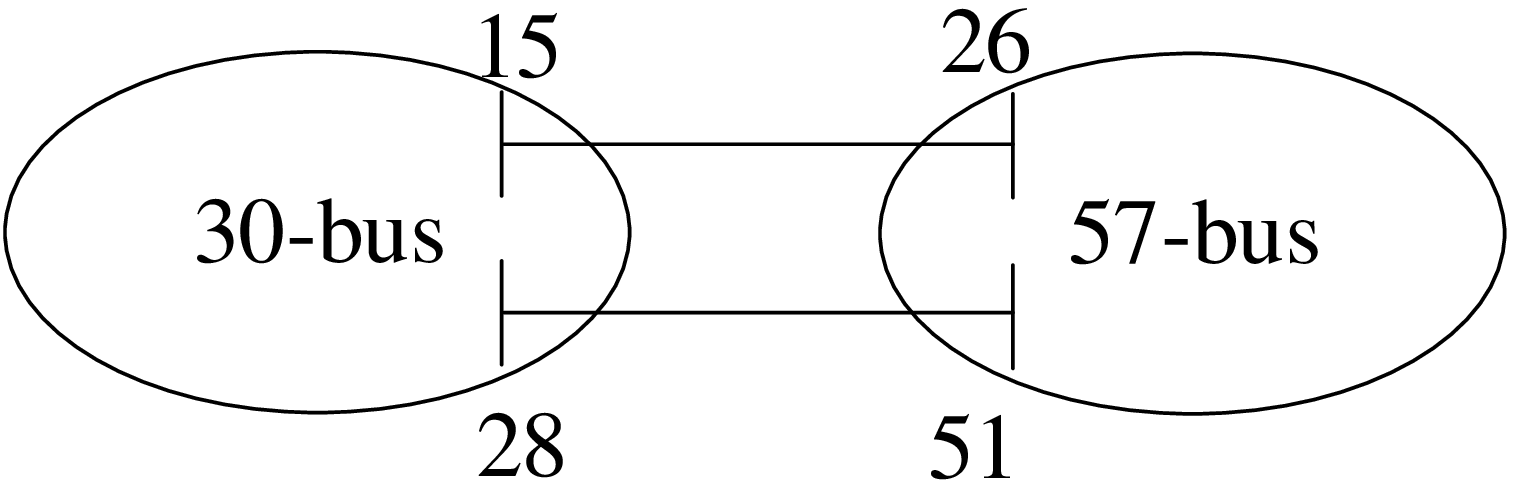}
        \caption{Two-area 87-bus system}
    \end{subfigure}
	\qquad
    \begin{subfigure}[t]{0.35\textwidth}
        \centering
        \includegraphics[width=\textwidth]{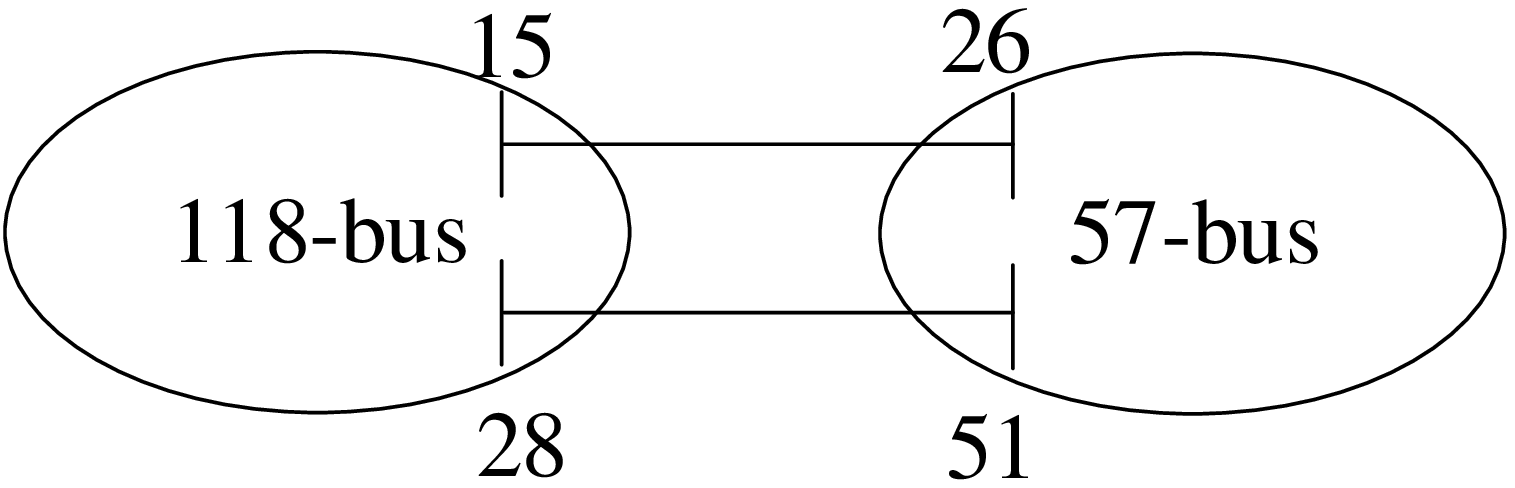}
        \caption{Two-area 175-bus system}
    \end{subfigure}

    \vspace{0.5cm}

    \begin{subfigure}[t]{0.35\textwidth}
        \centering
        \includegraphics[width=\textwidth]{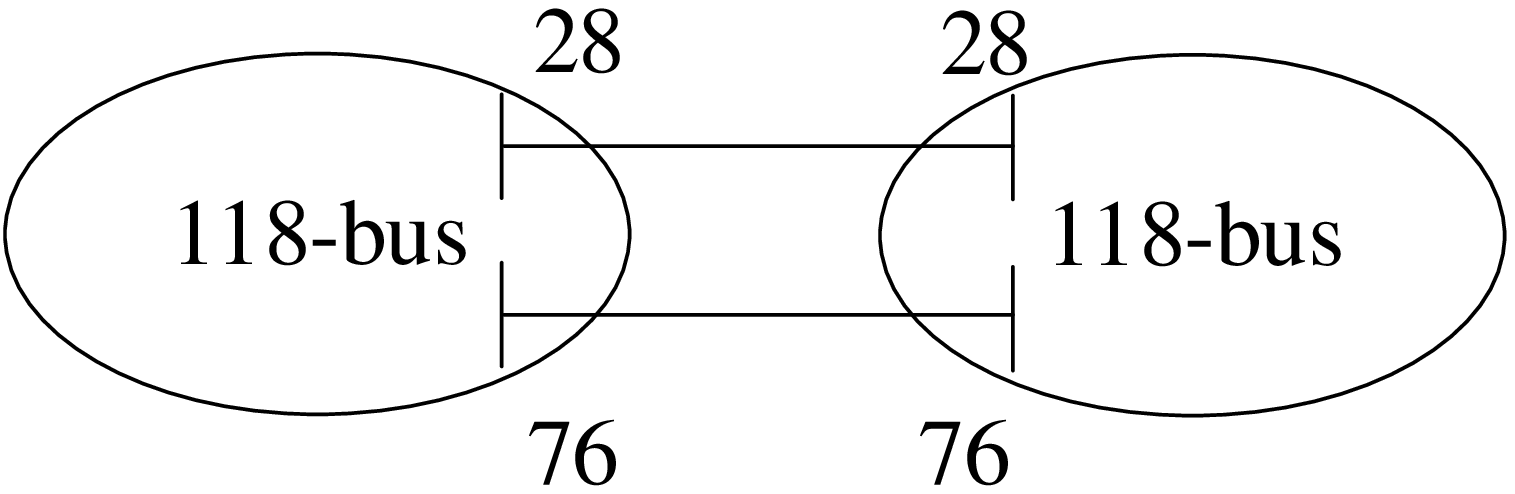}
        \caption{Two-area 236-bus system}
    \end{subfigure}
        \qquad
    \begin{subfigure}[t]{0.35\textwidth}
        \centering
        \includegraphics[width=\textwidth]{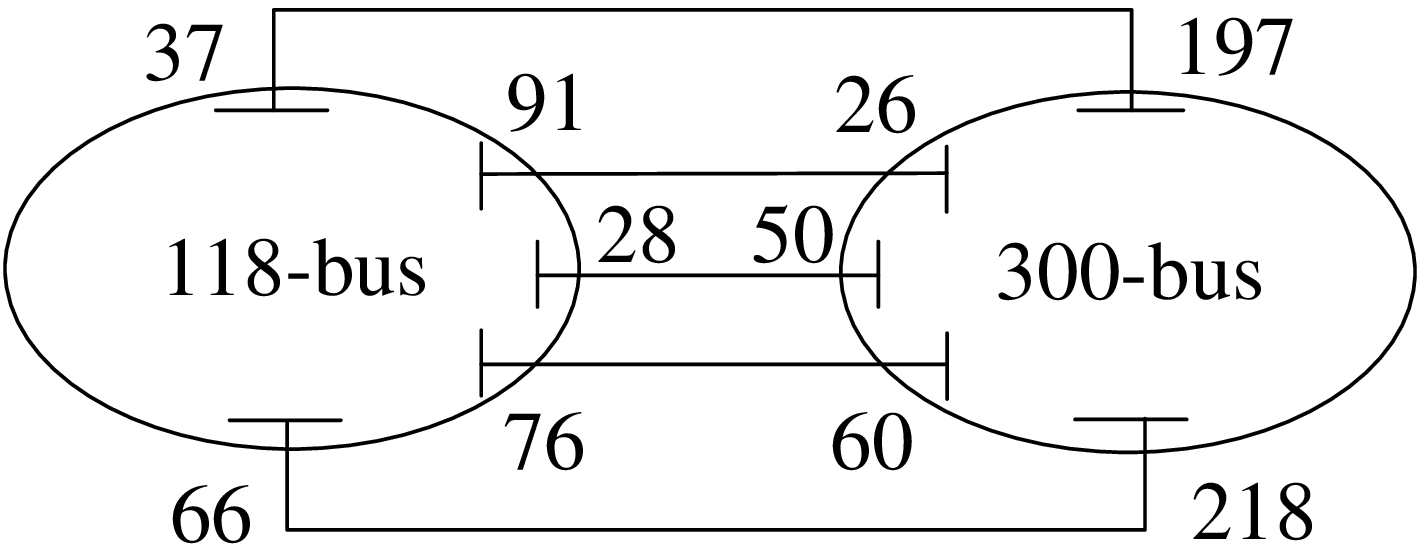}
        \caption{Two-area 418-bus system}
    \end{subfigure}

    \vspace{0.5cm}

    \begin{subfigure}[t]{0.35\textwidth}
        \centering
        \includegraphics[width=\textwidth]{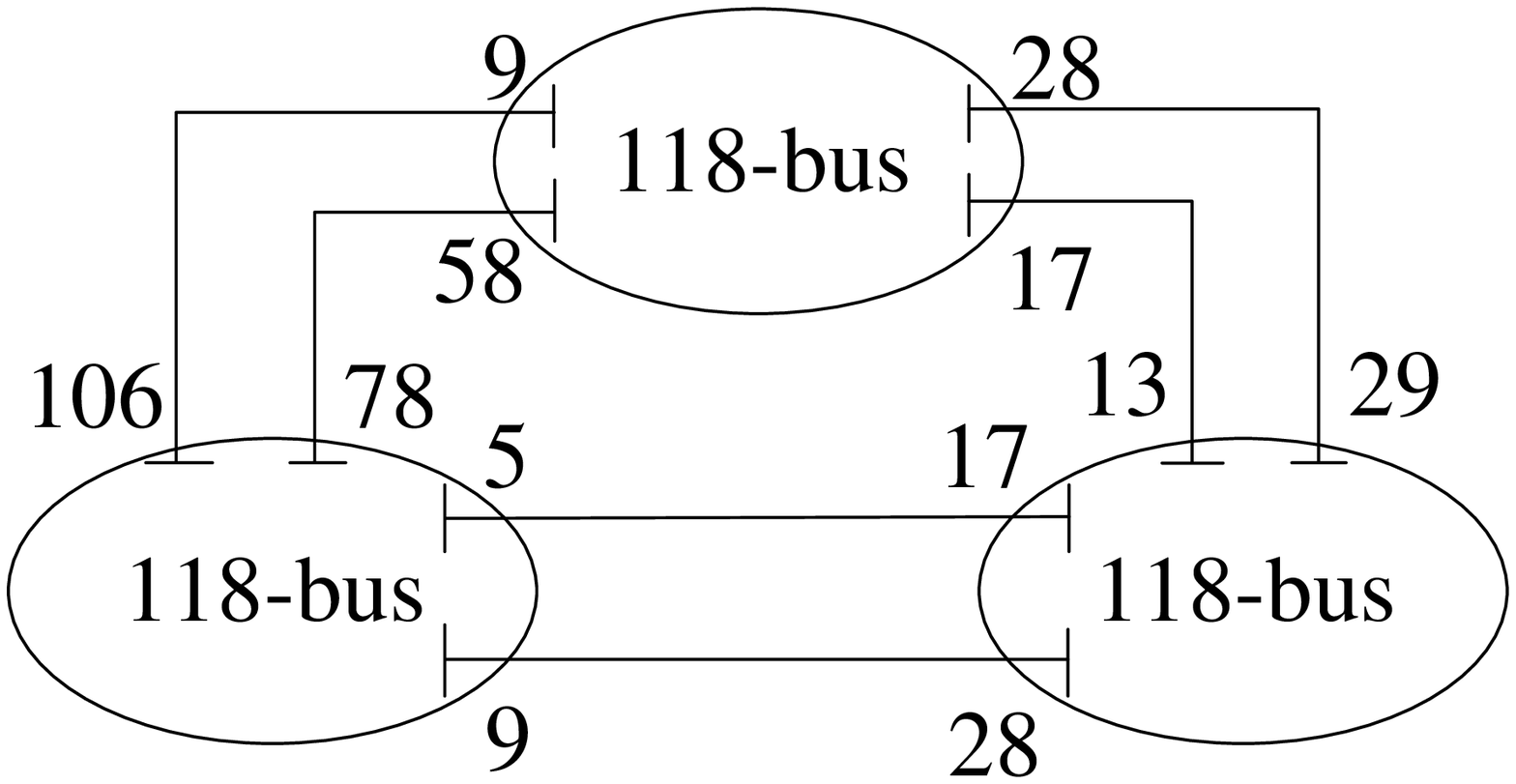}
        \caption{Three-area 354-bus system}
    \end{subfigure}
        \qquad
    \begin{subfigure}[t]{0.35\textwidth}
        \centering
        \includegraphics[width=\textwidth]{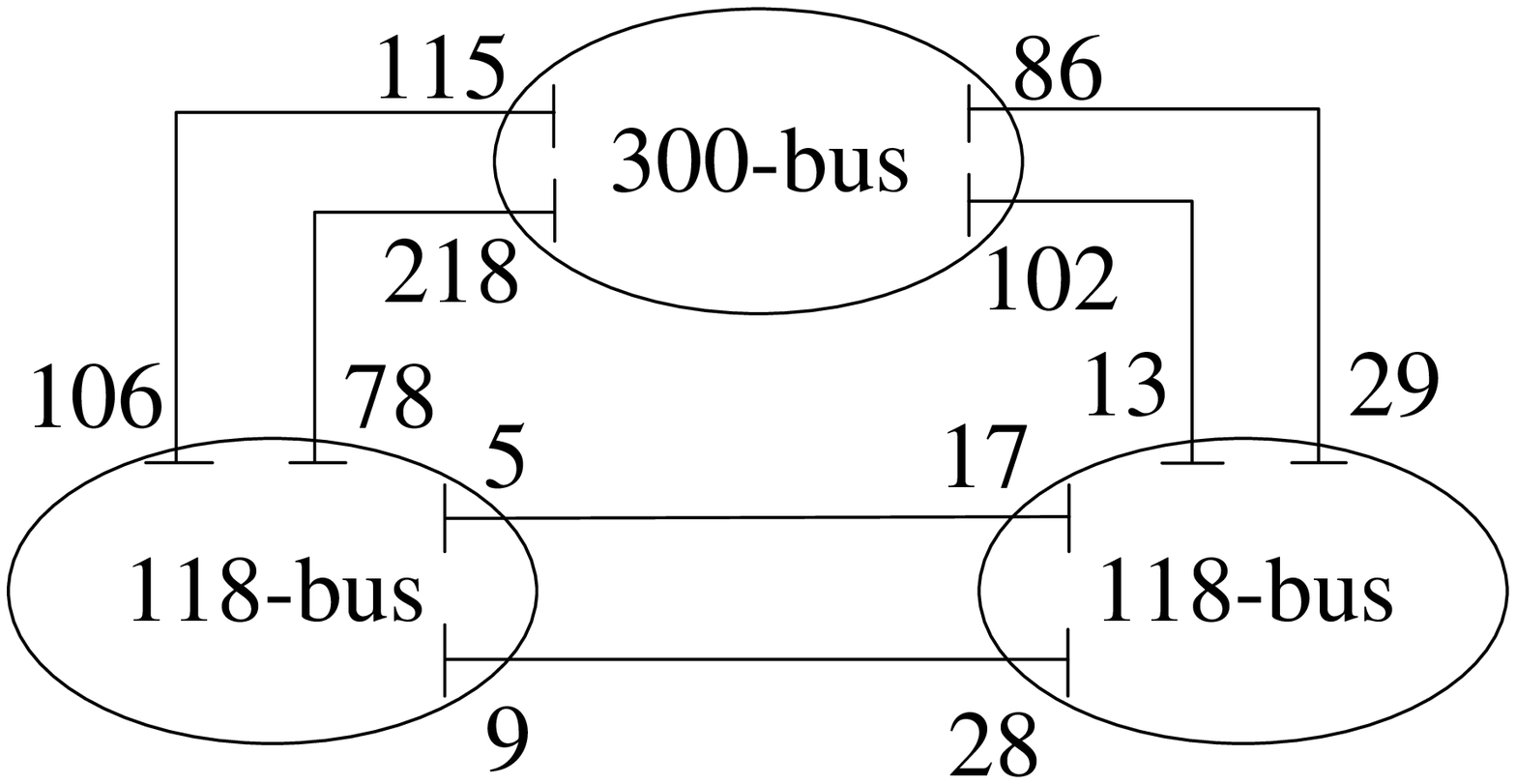}
        \caption{Three-area 536-bus system}
    \end{subfigure}
    \caption{Additional power system examples considered for numerical experiments.}
    \label{sec:additionalExamples}
\end{figure}

The multi-area power systems considered in Section \ref{sec:othertest} are given in Figure \ref{sec:additionalExamples}. Tie-line capacities were set to 100MW and their reactances were set to $0.25p.u.$ Capacity limits on the transmission lines within each area were set to their respective nominal values in Matpower \cite{MATPOWER} wherever present, and to 100MW, otherwise. For all two-area tests, two wind generators were installed in the two areas at buses 6 and 14 in area 1 and buses 11 and 23 in area 2. For the three-area tests, we replicated the placements described in Section \ref{sec:3area}. Power demands and available wind generations were varied the same way as in Sections \ref{sec:2area} and \ref{sec:3area}.

%
%
%
%
%
%
%
%
%
%
%
%
%

\end{document}